\documentclass[times]{iapress}
\usepackage{moreverb}

\usepackage[colorlinks,bookmarksopen,bookmarksnumbered,citecolor=red,urlcolor=red]{hyperref} % if run by PDFLaTex

\def\volumeyear{2019}
\def\volumenumber{7}
\def\volumemonth{March}
\usepackage{amsmath,amssymb}
\usepackage{graphicx}

\newtheorem{thm}{Theorem}[section]
\newtheorem{lem}{Lemma}[section]
\newtheorem{nas}{Corollary}[section]
\newtheorem{ozn}{Definition}[section]
\newtheorem{exm}{Example}[section]

\newcommand{\lp}[1]{\left( \begin{array}{#1} }
\newcommand{\rp}{\end{array} \right)}

\newcommand{\be}{\begin{equation}}
\newcommand{\ee}{\end{equation}}
\begin{document}

\runningheads{Extrapolation  Problem  for Sequences with Missing Observations}{ O. Masyutka, M. Moklyachuk and M. Sidei}

\title{Extrapolation   Problem  for Multidimensional Stationary Sequences with Missing Observations}

\author{Oleksandr Masyutka\affil{1}, Mikhail Moklyachuk \affil{2}$^,$\corrauth,
Maria Sidei\affil{2}}

\address{
\affilnum{1} Department of Mathematics and Theoretical Radiophysics,
Kyiv National Taras Shevchenko University, Kyiv 01601, Ukraine
\affilnum{2}Department of Probability Theory, Statistics and Actuarial
Mathematics, Taras Shevchenko National University of Kyiv, Kyiv 01601, Ukraine
}

\corraddr{Mikhail Moklyachuk (Email: Moklyachuk@gmail.com). Department of Probability Theory, Statistics and Actuarial
Mathematics, Taras Shevchenko National University of Kyiv, Volodymyrska 64 Str., Kyiv 01601, Ukraine.}

 \begin{abstract}
   This paper focuses on the  problem  of the mean square optimal estimation of linear functionals  which depend on the unknown values of a multidimensional stationary stochastic sequence.
   Estimates are based on observations of the sequence with an additive stationary noise sequence.
   The aim of the paper is to develop methods of finding the optimal estimates of the functionals in the case of missing observations.
   The problem is investigated in the case of spectral certainty where the spectral densities of the sequences are exactly known.
   Formulas for calculating the mean-square errors and the spectral characteristics of the optimal linear estimates of functionals are derived under the condition of spectral certainty.
   The minimax (robust) method of estimation is applied in the case of spectral uncertainty, where spectral densities of the sequences are not known exactly while sets of admissible spectral densities are given. Formulas that determine the least favorable spectral densities and the minimax spectral characteristics of the optimal estimates of functionals are proposed for some special sets of admissible densities.
\end{abstract}


\maketitle

\vspace{10pt}

\noindent {\bf DOI:} 10.19139/soic.v7i1.527

%--------------------------------------------------------------------------------

\section{Introduction}
The problem of estimation of the unknown values of stochastic processes is of constant interest in the theory and applications of stochastic processes.
The formulation of the estimation problems (interpolation, extrapolation and filtering) for stationary stochastic sequences with known spectral densities and reducing these problems to the corresponding problems of the theory of functions belongs to Kolmogorov \cite{Kolmogorov}. Effective methods of solution of the estimation problems for stationary stochastic sequences and processes were developed by Wiener \cite{Wiener} and Yaglom \cite{Yaglom1,Yaglom2}.  Further results are described in the books  by Rozanov \cite{Rozanov},  Hannan \cite{Hannan},
 Box et. al \cite{Box:Jenkins},  Brockwell and  Davis \cite{Brockwell_Davis}.
The crucial assumption of most of the methods developed for estimating of the unobserved values of stochastic processes is that the spectral densities of the involved stochastic processes are exactly known. In practice, however, complete information on the spectral densities is impossible in most cases.
In this situation one finds parametric or nonparametric estimates of the unknown spectral densities and then apply one of the traditional estimation methods provided that the selected spectral densities are true. This procedure can result in significant increasing of the value of the error of estimate as Vastola and Poor \cite{Vastola} have demonstrated with the help of some examples.
To avoid this effect one can search estimates which are optimal for all densities from a certain given class of admissible spectral densities. These estimates are called minimax since they minimize the maximum value of the error of estimates.
The paper by Grenander \cite{Grenander} was the first one where this approach to extrapolation problem for stationary processes was proposed.
Several models of spectral uncertainty and minimax-robust methods of data processing can be found in the survey paper by Kassam and Poor \cite{Kassam}.  Franke \cite{Franke1,Franke},  Franke and Poor \cite{Franke_Poor} investigated the minimax extrapolation  and filtering problems for stationary sequences with the help of convex optimization methods. This approach makes it possible to find equations that determine the least favorable spectral densities for some classes of admissible densities.
In the papers by Moklyachuk \cite{Moklyachuk:2008,Moklyachuk:2015} results of investigation of the  extrapolation, interpolation and filtering  problems for functionals which depend on the unknown values of stationary processes and sequences are described.
The problem of estimation of functionals which depend on the unknown values of multivariate stationary stochastic processes is the aim of the papers by
Moklyachuk and Masyutka \cite{Moklyachuk:Mas2008} - \cite{Moklyachuk:Mas2012}.
In the book  by Moklyachuk and Golichenko \cite{Golichenko} results of investigation of the interpolation, extrapolation and filtering problems for periodically correlated
stochastic sequences are proposed.
In their papers Luz and Moklyachuk \cite{Luz8} - \cite{Luz2016}  deal with the problems of estimation of functionals which depend on the unknown values of stochastic sequences with stationary increments.
Prediction problem for stationary sequences with missing observations is investigated in papers by Bondon \cite{Bondon1, Bondon2},
Cheng, Miamee and Pourahmadi \cite{Cheng1},
Cheng and Pourahmadi \cite{Cheng2},
Kasahara,  Pourahmadi and Inoue \cite{Kasahara},
Pourahmadi, Inoue and Kasahara \cite{Pourah},
Pelagatti \cite{Pelagatti}.
In papers by Moklyachuk and Sidei \cite{Sidei1} - \cite{Sidei7} an approach is developed to investigation of the interpolation, extrapolation and filtering  problems for stationary stochastic sequences with missing observations.

In this paper we present results of investigation of  the problem of the mean-square optimal estimation of the linear functional  $$A\vec{\xi}=\sum\limits_{j=0}^{\infty}\vec{a}(j)^\top\vec{\xi}(j)$$ which depends on the unknown values of a multivariate stationary stochastic sequence $\{\vec{\xi}(j),j\in\mathbb{Z}\}$.
Estimates are based on observations of the sequence with an additive stationary stochastic noise sequence $\{\vec{\xi}(j)+\vec{\eta}(j)\}$ at points $j\in\mathbb{Z}_{-}\backslash S=\{\ldots, -2,-1\}\backslash S$,
where $S=\bigcup\limits_{l=1}^{s}\{ -M_{l}-N_{l}, -M_{l}-N_{l}+1, \ldots,  -M_{l} \}$.
The problem is investigated in the case of spectral certainty, where the spectral densities of the signal and the noise  sequences $\{\vec{\xi}(j), j\in \mathbb{Z}\}$ and $\{\vec{\eta}(j), j\in  \mathbb{Z}\}$ are exactly known, and
in the case of spectral uncertainty, where the spectral densities of the sequences are not exactly known while a set of admissible spectral densities is given.
We first propose results of investigation of the mean-square optimal linear estimate of the linear functional in the case of spectral certainty. To find the optimal solution of the estimation problem in this case we apply an approach based on the Hilbert space projection method proposed by Kolmogorov \cite{Kolmogorov} and developed in the papers by Moklyachuk \cite{Moklyachuk:2008,Moklyachuk:2015}, and Moklyachuk and Masytka \cite{Moklyachuk:Mas2008} - \cite{Moklyachuk:Mas2012}.
We derive formulas for calculation the spectral characteristic and the mean-square error of the optimal estimate of the functional.
Next, in the case of spectral uncertainty, where the full information on spectral densities is impossible, while it is known that spectral densities of the sequences belong to some specified classes of admissible densities, the minimax-robust method of estimation is applied. This method gives us a procedure of finding estimates which minimize the maximum values of the mean-square errors of the estimates for all spectral densities from a given class of admissible spectral densities.
Formulas that determine the least favorable spectral densities and the minimax-robust spectral characteristics of the optimal estimates of the functional are proposed for some specific classes of admissible spectral densities.

\section{Hilbert space projection method of extrapolation of stationary sequences with missing observations}

Let $\vec{ \xi}(j)=\left \{ \xi_ {k} (j) \right \}_{k = 1} ^ {T},\,j\in \mathbb{Z}$ and $\vec{ \eta}(j)=\left \{ \eta_ {k} (j) \right \}_{k = 1} ^ {T},\,j\in \mathbb{Z}$, be multidimensional stationary stochastic sequences with zero mean values: $E\vec{\xi}(j)=\vec0$, $E\vec{\eta}(j)=\vec0$ and correlation functions which admit the spectral decomposition (see Gikhman and Skorokhod \cite{Gihman})
$$
R_{\xi}(n)=E\vec{\xi}(j+n)(\vec{\xi}(j))^{*}=\frac{1}{2\pi}\int\limits_{-\pi}^{\pi}e^{in\lambda}F(\lambda)d\lambda, \quad
  R_{\xi\eta}(n)=E\vec{\xi}(j+n)(\vec{\eta}(j))^{*}=\frac{1}{2\pi}\int\limits_{-\pi}^{\pi}e^{in\lambda}F_{\xi\eta}(\lambda)d\lambda,
 $$
  $$ R_{\eta\xi}(n)=E\vec{\eta}(j+n)(\vec{\xi}(j))^{*}=\frac{1}{2\pi}\int\limits_{-\pi}^{\pi}e^{in\lambda}F_{\eta\xi}(\lambda)d\lambda, \quad
  R_{\eta}(n)=E\vec{\eta}(j+n)(\vec{\eta}(j))^{*}=\frac{1}{2\pi}\int\limits_{-\pi}^{\pi}e^{in\lambda}G(\lambda)d\lambda,
$$
where $F(\lambda)=\left\{f_{kl}(\lambda)\right\}_{k,l=1}^T$, $F_{\xi\eta}(\lambda)=\left\{f_{kl}^{\xi\eta}(\lambda)\right\}_{k,l=1}^T$, $F_{\eta\xi}(\lambda)=\left\{f_{kl}^{\eta\xi}(\lambda)\right\}_{k,l=1}^T$, $G(\lambda)=\left\{g_{kl}(\lambda)\right\}_{k,l=1}^T$ are spectral densities of the stationary sequences such that the minimality condition holds true
\begin{equation}\label{minimal}
\int\limits_{-\pi}^{\pi}\text{Tr}\,\left(F(\lambda)+F_{\xi\eta}(\lambda)+F_{\eta\xi}(\lambda)+G(\lambda)\right)^{-1}d\lambda<\infty.
\end{equation}
Under this condition the mean-square error of the optimal estimate of the functional is nonzero (see Rozanov \cite{Rozanov}).

The stationary sequences  $\vec{\xi}(j)$ and $\vec{\eta}(j)$ admit the spectral decompositions  (see Gikhman and Skorokhod \cite{Gihman}; Karhunen \cite{Karhunen})
\begin{equation} \label{ksi}
\vec{\xi}(j)=\int\limits_{-\pi}^{\pi}e^{ij\lambda}Z_{\xi}(d\lambda), \hspace{1cm}
\vec{\eta}(j)=\int\limits_{-\pi}^{\pi}e^{ij\lambda}Z_{\eta}(d\lambda),
\end{equation}
where $Z_{\xi}(d\lambda)$ and $Z_{\eta}(d\lambda)$ are orthogonal stochastic measures  such that the following relations hold true
$$
EZ_{\xi}(\Delta_1)(Z_{\xi}(\Delta_2))^*=\frac{1}{2\pi}\int_{\Delta_1\cap\Delta_2}F(\lambda)d\lambda, \quad
EZ_{\xi}(\Delta_1)(Z_{\eta}(\Delta_2))^*=\frac{1}{2\pi}\int_{\Delta_1\cap\Delta_2}F_{\xi\eta}(\lambda)d\lambda,
$$
$$
EZ_{\eta}(\Delta_1)(Z_{\xi}(\Delta_2))^*=\frac{1}{2\pi}\int_{\Delta_1\cap\Delta_2}F_{\eta\xi}(\lambda)d\lambda, \quad
EZ_{\eta}(\Delta_1)(Z_{\eta}(\Delta_2))^*=\frac{1}{2\pi}\int_{\Delta_1\cap\Delta_2}G(\lambda)d\lambda.
$$

Consider the problem of the mean-square optimal linear estimation of the functional
$$A\vec{\xi}=\sum\limits_{j=0}^{\infty}\vec a(j)^\top\vec{\xi}(j)$$
which depends on the unknown values of the sequence $\{\vec{\xi}(j),j\in \mathbb{Z}\}$ from observations of the sequence $\vec{\xi}(j)+\vec{\eta}(j)$ at points $j\in\mathbb{Z}_{-}\backslash S $, where $S=\bigcup\limits_{l=1}^{s}\{ -M_{l}-N_{l}, \ldots, -M_{l} \}.$

Making use of the spectral decomposition (\ref{ksi}) of the sequence $\vec{\xi}(j)$ we can represent the functional $A\vec{\xi}$ in the form
\begin{equation*}
A\vec{\xi}=\int\limits_{-\pi}^{\pi}(A(e^{i\lambda}))^\top Z_{\xi}(d\lambda), \quad
  A(e^{i\lambda})=\sum\limits_{j=0}^{\infty}\vec a(j)e^{ij\lambda}.
 \end{equation*}

We will suppose that the coefficients $\{\vec a(j), j= 0, 1, \ldots\}$ which determine the functional  $A\vec{\xi}$ are such that the following condition
\begin{equation}\label{umovu2mom}
\sum\limits_{j=0}^{\infty}\sum_{k=1}^T\left|a_k(j)\right|<\infty
\end{equation}
is satisfied. This condition ensures that the functional  $A\vec{\xi}$ has a finite second moment.

Denote by $\hat{A}\vec{\xi}$ the optimal linear estimate of the functional  $A\vec{\xi}$ from the known observations of the sequence  $\vec{\xi}(j)+\vec{\eta}(j)$ at points $j\in\mathbb{Z}_{-}\backslash S $.
 Since the spectral densities of the stationary sequences   $\vec{\xi}(j)$ and $\vec{\eta}(j)$ are suppose to be known, we can use the Hilbert space projection method proposed by A.~N.~Kolmogorov (see selected works by Kolmogorov \cite{Kolmogorov}) to find the estimate  $\hat{A}\vec{\xi}$.

Consider  values $\xi_k(j),k= {1,\dots,T},j\in\mathbb Z$, and $\eta_k(j),k= {1,\dots,T},j\in\mathbb Z$, of the sequences as  elements of the Hilbert space $H=L_2(\Omega,\mathcal{F},P)$ generated by random variables $\xi$ with zero mathematical expectations, $E\xi=0$,  finite variations, $E|\xi|^2<\infty$, and the inner product $(\xi,\eta)=E(\xi\overline{\eta})$.
Denote by $H^s(\xi+\eta)$ the closed linear subspace generated by elements $\{\xi_k(j)+\eta_k(j): j\in \mathbb{Z}_{-}\backslash S, k=\overline{1,T}\}$ in the Hilbert space $H=L_2(\Omega,\mathcal{F},P)$.
Denote by $L_{2} (F+G)$ the Hilbert space of vector-valued functions $ \vec{a}( \lambda )= \left \{a_{k} ( \lambda ) \right \}_{k=1}^{T} $ such that
 \[ \int_{- \pi}^{ \pi} \vec{a}( \lambda )^{ \top} \left(F(\lambda)+F_{\xi\eta}(\lambda)+F_{\eta\xi}(\lambda)+G(\lambda)\right) \overline{ \vec{a}( \lambda )}d \lambda < \infty. \] Denote by $L_2^s(F+G)$ the subspace of $L_2(F+G)$ generated by functions of the form
 \[e^{in \lambda} \delta_{k} , \; \delta_{k} = \left \{ \delta_{kl} \right \}_{l=1}^{T} , \; k= {1,\dots,T}, \; n \in Z_{-} \backslash S. \]

The mean-square optimal linear estimate $\hat{A}\vec{\xi}$ of the functional $A\vec{\xi}$ from observations of the sequence $\vec{\xi}(j)+\vec{\eta}(j)$ is of the form
 \begin {equation}\label{sp-har-general}
\hat{A}\vec{\xi}=\int\limits_{-\pi}^{\pi}(h(e^{i\lambda}))^\top(Z_{\xi}(d\lambda)+ Z_{\eta}(d\lambda)),
 \end{equation}
where $h(e^{i\lambda})=\left\{h_k(e^{i\lambda})\right\}_{k=1}^T\in L_2^s(F+G)$ is the spectral characteristic of the estimate.

The mean-square error  of the estimate $\hat{A}\vec{\xi}$ is given by the formula
\[\Delta(h;F,G)=E\left|A\vec{\xi}-\hat{A}\vec{\xi}\right|^2=\]
 \[=\frac{1}{2\pi}\int\limits_{-\pi}^{\pi}\left(A(e^{i\lambda})-h(e^{i\lambda})\right)^\top F(\lambda)\overline{\left(A(e^{i\lambda})-h(e^{i\lambda})\right)}d\lambda
 +\frac{1}{2\pi}\int\limits_{-\pi}^{\pi}\left(h(e^{i\lambda})\right)^\top G(\lambda)\overline{\left(h(e^{i\lambda})\right)}d\lambda-\]
\begin{equation}\label{error-general}
-\frac{1}{2\pi}\int\limits_{-\pi}^{\pi}\left(A(e^{i\lambda})-h(e^{i\lambda})\right)^\top F_{\xi\eta}(\lambda)\overline{h(e^{i\lambda})}d\lambda -\frac{1}{2\pi}\int\limits_{-\pi}^{\pi}\left(h(e^{i\lambda})\right)^\top F_{\eta\xi}(\lambda)\overline{\left(A(e^{i\lambda})-h(e^{i\lambda})\right)}d\lambda
\end{equation}

According to the Hilbert space orthogonal projection method the optimal linear estimate of the functional $A\vec{\xi}$  is a projection of the element $A\vec{\xi}$ of the space $H$  on  the subspace $H^s(\xi+\eta)$. The projection is determined by  the following conditions:
\begin{equation*} \begin{split}
1)& \hat{A}\vec{\xi} \in H^s(\xi+\eta), \\
2)& A\vec{\xi}-\hat{A}\vec{\xi} \bot  H^s(\xi+\eta).
\end{split} \end{equation*}

It follows from the second condition that the spectral characteristic $h(e^{i\lambda})=\left\{h_k(e^{i\lambda})\right\}_{k=1}^T$ of the optimal linear estimate $\hat{A}\vec{\xi}$ for any $ j\in \mathbb{Z}_{-}\backslash S $ satisfies  equations
\begin{equation*} \begin{split}
&\frac{1}{2\pi}\int\limits_{-\pi}^{\pi} \left(A(e^{i\lambda})- h(e^{i\lambda})\right)^\top F(\lambda)e^{-ij\lambda}d\lambda -\frac{1}{2\pi}\int\limits_{-\pi}^{\pi} (h(e^{i\lambda}))^\top {F_{\eta\xi}}(\lambda)e^{-ij\lambda}d\lambda +\\
&+\frac{1}{2\pi}\int\limits_{-\pi}^{\pi} \left(A(e^{i\lambda})- h(e^{i\lambda})\right)^\top F_{\xi\eta}(\lambda)e^{-ij\lambda}d\lambda -\frac{1}{2\pi}\int\limits_{-\pi}^{\pi} (h(e^{i\lambda}))^\top G(\lambda)e^{-ij\lambda}d\lambda=\vec{0}.
\end{split} \end{equation*}

The last relation can be written in the form
\begin{equation*} \begin{split}
\frac{1}{2\pi}\int\limits_{-\pi}^{\pi} \left[(A(e^{i\lambda}))^\top(F(\lambda)+F_{\xi\eta}(\lambda))-(h(e^{i\lambda}))^\top(F(\lambda)+F_{\xi\eta}(\lambda)+F_{\eta\xi}(\lambda) +G(\lambda))\right]e^{-ij\lambda}d\lambda=\vec{0},\; j\in \mathbb{Z}_{-}\backslash S.
\end{split} \end{equation*}

Hence the function $\left[(A(e^{i\lambda}))^\top(F(\lambda)+F_{\xi\eta}(\lambda))-(h(e^{i\lambda}))^\top(F(\lambda)+F_{\xi\eta}(\lambda)+F_{\eta\xi}(\lambda) +G(\lambda))\right]$ is of the form

$$(A(e^{i\lambda}))^\top(F(\lambda)+F_{\xi\eta}(\lambda))-(h(e^{i\lambda}))^\top(F(\lambda)+F_{\xi\eta}(\lambda)+F_{\eta\xi}(\lambda) +G(\lambda))=(C(e^{i\lambda}))^\top,$$
$$C(e^{i\lambda})=\sum\limits_{j\in U} \vec{c}(j)e^{ij\lambda},$$
where $U=S\cup \{0,1, \ldots\},$ and  $\vec{c}(j), j \in U$ are unknown coefficients to be determined.

From the last relation we deduce that the spectral characteristic of the optimal linear estimate $\hat{A}\vec{\xi}$ is of the form
\begin{equation} \label{sphar} \begin{split}
(h(e^{i\lambda}))^\top=(A(e^{i\lambda}))^\top(F(\lambda)+ F_{\xi\eta}(\lambda))(F_{\zeta}(\lambda))^{-1}-(C(e^{i\lambda}))^\top(F_{\zeta}(\lambda))^{-1},
\end{split} \end{equation}
where $F_{\zeta}(\lambda)=F(\lambda)+F_{\xi\eta}(\lambda)+F_{\eta\xi}(\lambda)+G(\lambda).$

From the first condition, $\hat{A}\vec{\xi}\in H^s(\xi+\eta)$, which determines the optimal estimate of the functional  $A\vec{\xi}$, it follows that
$$\frac{1}{2\pi}\int\limits_{-\pi}^{\pi} h(e^{i\lambda})e^{-ij\lambda }d\lambda=\vec{0}, \, j \in U,$$
namely
\begin{equation*} \begin{split}
\frac{1}{2\pi}\int\limits_{-\pi}^{\pi}\left(
(A(e^{i\lambda}))^\top(F(\lambda)+ F_{\xi\eta}(\lambda))(F_{\zeta}(\lambda))^{-1}-(C(e^{i\lambda}))^\top(F_{\zeta}(\lambda))^{-1}\right)e^{-ij\lambda}d\lambda=\vec{0},\; j\in  U.
\end{split} \end{equation*}

Disclose brackets and write the last equation in the form
\begin{equation}\label{3} \begin{split}
\sum\limits_{k=0}^{\infty}(\vec{a}(k))^\top\frac{1}{2\pi}\int\limits_{-\pi}^{\pi}(F(\lambda)+ F_{\xi\eta}(\lambda))(F_{\zeta}(\lambda))^{-1}e^{i(k-j)\lambda}d\lambda
-\sum\limits_{l\in U}(\vec{c}(l))^\top\frac{1}{2\pi}\int\limits_{-\pi}^{\pi}(F_{\zeta}(\lambda))^{-1}e^{i(l-j)\lambda}d\lambda=\vec{0}.
\end{split} \end{equation}

Let us introduce the Fourier coefficients of the functions

\begin{equation}\begin{split}\label{Fcoef}
&B(k-j)=\frac{1}{2\pi} \int\limits_{-\pi}^{\pi}(F_{\zeta}(\lambda))^{-1}e^{-i(k-j)\lambda}d\lambda;\\
&R(k-j)=\frac{1}{2\pi} \int\limits_{-\pi}^{\pi}(F(\lambda)+ F_{\xi\eta}(\lambda))(F_{\zeta}(\lambda))^{-1}e^{-i(k-j)\lambda}d\lambda;\\
&Q(k-j)=\frac{1}{2\pi} \int\limits_{-\pi}^{\pi}F(\lambda)(F_{\zeta}(\lambda))^{-1}G(\lambda)-F_{\xi\eta}(\lambda)(F_{\zeta}(\lambda))^{-1}F_{\eta\xi}(\lambda)e^{-i(k-j)\lambda}d\lambda.\\
\end{split}\end{equation}

Denote by
$\vec{\bold{a}}^\top=(0, 0, \ldots, 0, \vec{a}^\top)$ a vector that has first $T\cdot|S|=T\cdot\sum\limits_{k=1}^{s}(N_k+1)$ zero components, and the last component $\vec{a}^\top=(\vec{a}(0)^\top,\vec{a}(1)^\top\,\ldots)$ is constructed of  coefficients which define the functional  $A\vec{\xi}$.

Now we can represent relation (\ref{3}) in the form
\begin{equation}\label{rivn2}\begin{split}
\bold{R}\vec{\bold{a}}=\bold{B} \vec{\bold{c}},
\end{split} \end{equation}
where $\vec{\bold{c}}$ is the vector constructed of the unknown coefficients $\vec{c}(k),k\in U$. The linear operator $\bold{B}$ in the space  $\ell_2$ is defined by the matrix
  \begin{equation*} \label{matrix2}
B= \left( \begin{array}{ccccc}
B_{s, s}& B_{s, s-1}&\ldots&B_{s, 1}&B_{s, n}\\
B_{s-1, s}& B_{s-1, s-1}&\ldots&B_{s-1, 1}&B_{s-1, n}\\
\vdots&\vdots&\ddots &\vdots&\vdots\\
B_{1, s}& B_{1, s-1}&\ldots&B_{1, 1}&B_{1, n}\\
B_{n, s}& B_{n, s-1}&\ldots&B_{n, 1}&B_{n, n}
\end{array}\right),
\end{equation*}
 where   elements in the last column and the last row are compound matrices constructed with the help of the block-matrices
  \begin{equation*}\begin{split}
&B_{l,n}(k,j)= B(k-j),\; l= 1,2,\ldots,s;\; k= -M_{l}-N_l, \ldots,-M_{l};\;  j= 0, 1, 2, \ldots,\\
&B_{n,m}(k,j)=B(k-j),\; m= 1,2,\ldots,s;\; k= 0, 1, 2,\ldots;\; j= -M_{m}-N_m, \ldots,-M_{m},\\
& B_{n,n}(k,j)=B(k-j),\; k,j= 0, 1, 2, \ldots,
\end{split} \end{equation*}
and other elements of  matrix $B$ are compound matrices constructed with the help of the block-matrices
\begin{equation*} \label{matr32}\begin{split}
&B_{l,m}(j,k)=B(k-j),\; l,m=1,2,\ldots,s;\; k= -M_{l}-N_l, \ldots,-M_{l};\;  j= -M_{m}-N_m, \ldots,-M_{m}.
\end{split}\end{equation*}

The linear operator $\bold{R}$  in the space  $\ell_2$ is defined by the corresponding matrix in the same manner.

 The unknown coefficients $\vec{c}(k),k\in U$, which are determined by equation (\ref{rivn2}) can be calculated by the formula
$$\vec{c}(k)=(\bold{B}^{-1}\bold{R}\vec{\bold{a}})(k),$$
where $(\bold{B}^{-1}\bold{R}\vec{\bold{a}})(k)$ is the $k$-th component  of the vector $\bold{B}^{-1}\bold{R}\vec{\bold{a}}$.
We will suppose that  the operator $\bold{B}$ is invertible (see paper by Salehi \cite{Salehi} for more details).

Hence the spectral characteristic   $h(e^{i\lambda})$ of the estimate  $\hat{A}\vec{\xi}$ can be calculated by the formula
\begin{equation}\label{4} \begin{split}
(h(e^{i\lambda}))^\top=(A(e^{i\lambda}))^\top(F(\lambda)+ F_{\xi\eta}(\lambda))(F_{\zeta}(\lambda))^{-1}-
\left(\sum\limits_{k\in T}(\bold{B}^{-1}\bold{R}\vec{\bold{a}})(k)e^{ik\lambda}\right)^\top(F_{\zeta}(\lambda))^{-1}.
\end{split} \end{equation}

The mean-square error of the estimate  $\hat{A}\xi$  can be calculated by the formula (\ref{error-general}) which can be represented in the form
\[\Delta(h;F,G,F_{\xi\eta},F_{\eta\xi})=E\left|A\vec{\xi}-\hat{A}\vec{\xi}\right|^2=\]
 \[=\frac{1}{2\pi}\int\limits_{-\pi}^{\pi}(a(\lambda))^\top F(\lambda)\overline{a(\lambda)}d\lambda+\frac{1}{2\pi}\int\limits_{-\pi}^{\pi}(b(\lambda))^\top G(\lambda)\overline{b(\lambda}d\lambda-\]
 \[-\frac{1}{2\pi}\int\limits_{-\pi}^{\pi}(a(\lambda))^\top F_{\xi\eta}(\lambda)\overline{b(\lambda}d\lambda-\frac{1}{2\pi}\int\limits_{-\pi}^{\pi}(b(\lambda))^\top F_{\eta\xi}(\lambda)\overline{a(\lambda)}d\lambda=\]
\[=\frac{1}{2\pi}\int\limits_{-\pi}^{\pi}(A(e^{i\lambda}))^\top(F(\lambda)G(\lambda)-F_{\xi\eta}(\lambda)F_{\eta\xi}(\lambda))(F_{\zeta}(\lambda))^{-1}\overline{A(e^{i\lambda})}d\lambda
+\frac{1}{2\pi}\int\limits_{-\pi}^{\pi}(C(e^{i\lambda}))^\top(F_{\zeta}(\lambda))^{-1}\overline{C(e^{i\lambda})}d\lambda=\]
\begin{equation}\label{55}
= \langle \bold{R}\vec{\bold{a}},\bold{B}^{-1}\bold{R}\vec{\bold{a}}\rangle +\langle {\bold{Q}\vec{\bold{a}},\vec{\bold{a}}}\rangle,
\end{equation}
where
\[(a(\lambda))^\top=(A(e^{i\lambda}))^\top(F_{\eta\xi}(\lambda)+G(\lambda))(F_{\zeta}(\lambda))^{-1}+(C(e^{i\lambda}))^\top(F_{\zeta}(\lambda))^{-1},\]
\[(b(\lambda))^\top=(A(e^{i\lambda}))^\top(F(\lambda)+ F_{\xi\eta}(\lambda))(F_{\zeta}(\lambda))^{-1}-(C(e^{i\lambda}))^\top(F_{\zeta}(\lambda))^{-1},\]
and $\langle {a,c}\rangle = \sum\limits_{k}a_k\overline{c}_k$\, is the inner product in the space $ \ell_2$.

The linear operator  $\bold{Q}$ in the space  $\ell_2$ is defined by the corresponding matrix in the same manner as operator  $\bold{B}$ is defined.

Thus we obtain the following theorem.

\begin{thm}\label{t2}
 Let $\{\vec{\xi}(j), j\in \mathbb{Z}\}$ and $\{\vec{\eta}(j), j\in  \mathbb{Z}\}$ be multidimensional stationary stochastic sequences with the spectral density matrices $F(\lambda), F_{\xi\eta}(\lambda), F_{\eta\xi}(\lambda), G(\lambda)$ and let the minimality condition (\ref{minimal}) be satisfied. The spectral characteristic   $h(e^{i\lambda})$ and the mean-square error  $\Delta(h;F,G,F_{\xi\eta},F_{\eta\xi})$ of the optimal linear estimate of the functional  $A\vec{\xi}$ which depends on the  unknown values of the sequence  $\vec{\xi}(j)$ based on observations of the sequence $\vec{\xi}(j)+\vec{\eta}(j),$ $j\in \mathbb{Z}_{-}\backslash S$ can be calculated by formulas (\ref{4}),  (\ref{55}).
\end{thm}

%------------------ corresponding results for uncorrelated sequences---------------------
The corresponding results can be obtained for the uncorrelated sequences $\{\vec{\xi}(j), j\in \mathbb{Z}\}$ and $\{\vec{\eta}(j), j\in  \mathbb{Z}\}$. In this case the spectral densities $F_{\xi\eta}(\lambda)=0$, $F_{\eta\xi}(\lambda)=0$ and we get the following corollary.

\begin{nas} \label{c21}
Let $\{\vec{\xi}(j), j\in \mathbb{Z}\}$ and $\{\vec{\eta}(j), j\in  \mathbb{Z}\}$ be uncorrelated multidimensional stationary stochastic sequences with spectral densities $F(\lambda)$ and $G(\lambda)$ which satisfy the minimality condition
\begin{equation}\label{minimal-Unc}
\int\limits_{-\pi}^{\pi}\text{Tr}\,(F(\lambda)+G(\lambda))^{-1}d\lambda<\infty.
\end{equation}
 The spectral characteristic   $h(e^{i\lambda})$ and the mean-square error  $\Delta(F,G)$ of the optimal linear estimate of the functional  $A\vec{\xi}$ which depends on  unknown values of the sequence  $\vec{\xi}(j)$ based on observations of the sequence $\vec{\xi}(j)+\vec{\eta}(j),$ $j\in \mathbb{Z}_{-}\backslash S$ can be calculated by the formulas
  \begin{equation}\label{sp-ch2} \begin{split}
(h(e^{i\lambda}))^\top=(A(e^{i\lambda}))^\top F(\lambda)(F(\lambda)+G(\lambda))^{-1}-\left(\sum\limits_{k\in U}(\bold{B}^{-1}\bold{R}\vec{\bold{a}})(k)e^{ik\lambda}\right)^\top(F(\lambda)+G(\lambda))^{-1},
\end{split} \end{equation}
\[\Delta(h;F,G)=E\left|A\vec{\xi}-\hat{A}\vec{\xi}\right|^2=\frac{1}{2\pi}\int\limits_{-\pi}^{\pi}(r_G(\lambda))^\top F(\lambda)\overline{r_G(\lambda)}d\lambda+\frac{1}{2\pi}\int\limits_{-\pi}^{\pi}(r_F(\lambda))^\top G(\lambda)\overline{r_F(\lambda)}d\lambda=\]
\begin{equation} \label{err2}
=\langle\bold{R}\vec{\bold{a}},\bold{B}^{-1}\bold{R}\vec{\bold{a}}\rangle+\langle\bold{Q}\vec{\bold{a}},\vec{\bold{a}}\rangle,
\end{equation}
where
\[(r_F(\lambda))^\top=\left((A(e^{i\lambda}))^\top F(\lambda)-\left(\sum\limits_{k \in U}(\bold{B}^{-1}\bold{R}\vec{\bold{a}})(k) e^{ik\lambda}\right)^\top\right)(F(\lambda)+G(\lambda))^{-1},\]
\[(r_G(\lambda))^\top=\left((A(e^{i\lambda}))^\top G(\lambda)+\left(\sum\limits_{k \in U}(\bold{B}^{-1}\bold{R}\vec{\bold{a}})(k) e^{ik\lambda}\right)^\top\right)(F(\lambda)+G(\lambda))^{-1},\]
and $\bold{B}, \bold{R}, \bold{Q}$ are linear operators in the space  $\ell_2$ that are determined by compound matrices constructed  of the block-matrices $B(k-j), R(k-j), Q(k-j)$ respectively which are  defined by  the Fourier coefficients of the functions

\begin{equation}\begin{split}\label{Fcoef-Un}
&B(k-j)=\frac{1}{2\pi} \int\limits_{-\pi}^{\pi}(F(\lambda)+G(\lambda))^{-1}e^{-i(k-j)\lambda}d\lambda;\\
&R(k-j)=\frac{1}{2\pi} \int\limits_{-\pi}^{\pi}F(\lambda)(F(\lambda)+G(\lambda))^{-1}e^{-i(k-j)\lambda}d\lambda;\\
&Q(k-j)=\frac{1}{2\pi} \int\limits_{-\pi}^{\pi}F(\lambda)(F(\lambda)+G(\lambda))^{-1}G(\lambda)e^{-i(k-j)\lambda}d\lambda.\\
\end{split}\end{equation}
\end{nas}

% -------------------------without noise------------

 Consider the estimation problem in the case where the stationary sequence  $\{\vec{\xi}(j), j\in \mathbb{Z}\}$ is observed without noise. Since in this case  $G(\lambda)=0$, the spectral characteristic of the estimate $\hat{A}\vec{\xi}$ is of the form
 \begin{equation} \label{sphar1}
(h(e^{i\lambda}))^\top=(A(e^{i\lambda}))^\top-(C(e^{i\lambda}))^\top (F(\lambda))^{-1},\quad C(e^{i\lambda})=\sum\limits_{j\in U} \vec{c}(j)e^{ij\lambda},
\end{equation}
and the system of equations (\ref{rivn2}) can  be represented in the form
\begin{equation} \label{rivn}
 \vec{\bold{a}}=\bold{B}\vec{\bold{c}},
 \end{equation}
where  $\bold{B}$ is the linear operator in the space $\ell_2$ which is constructed with the help of  the Fourier coefficients of the function $(F(\lambda))^{-1}$ and is of the similar form as operators defined before.

Hence, the unknown coefficients $\vec{c}(j), j \in U,$ can be calculated by the formula
 $$\vec{c}(j)=\left(\bold{B}^{-1}\vec{\bold{a}}\right)(j), $$
where $\left(\bold{B}^{-1}\vec{\bold{a}}\right)(j) $ is the $j$-th component of the vector $\bold{B}^{-1}\vec{\bold{a}},$ and the spectral characteristic of the estimate $\hat{A}\vec{\xi}$ is determined by the formula
\begin {equation} \label{spchar2}
(h(e^{i\lambda})^\top=\left(\sum\limits_{j=0}^{\infty}\vec{a}(j)e^{ij\lambda }\right)^\top-\left(\sum\limits_{j \in U}\left(\bold{B}^{-1}\vec{\bold{a}}\right)(j) e^{ij\lambda }\right)^\top(F(\lambda))^{-1}.
\end{equation}
The mean-square error of the estimate  $\hat{A}\vec{\xi}$ is determined by the formula
\begin {equation} \label {6} \begin{split}
 \Delta(h;F)=\langle\bold{B}^{-1}\vec{\bold{a}}, \vec{\bold{a}}\rangle.
\end{split} \end{equation}

Let us summarize the obtained result in the form of a corollary.

\begin{nas} \label{c22}
Let $\{\vec{\xi}(j), j\in \mathbb{Z}\}$ be a multidimensional stationary stochastic sequence with the spectral density $F(\lambda)$ which satisfy the  minimality condition
\begin{equation}\label{minimal1}
\int\limits_{-\pi}^{\pi}\text{Tr}\,(F(\lambda))^{-1}d\lambda<\infty.
\end{equation}
The spectral characteristic $h(e^{i\lambda})$ and the mean-square error   $\Delta(h,F)$ of the optimal linear estimate $\hat{A}\vec{\xi}$ of the functional $A\vec{\xi}$ from  observations of the sequence $\vec{\xi}(j)$ at points  $j\in\mathbb{Z}_{-}\backslash S $, where $S=\bigcup\limits_{l=1}^{s}\{-M_{l}-N_l, \ldots,-M_{l} \},$ can be calculated by formulas (\ref{spchar2}), (\ref{6}).
\end{nas}

%   -------------------------               A_N(\xi)     ----
Let $\vec{\xi}(j)$ and $\vec{\eta}(j)$ be uncorrelated stationary sequences.
Consider the problem of the mean-square optimal linear extrapolation of the functional
$$A_N\vec{\xi}=\sum\limits_{j=0}^{N}\vec{a}(j)^\top\vec{\xi}(j)$$
which depends on  unknown values of the sequence $\vec{\xi}(j)$  from  observations of the sequence  $\vec{\xi}(j)+\vec{\eta}(j)$ at points $j\in\mathbb{Z}_{-}\backslash S $, where $S=\bigcup\limits_{l=1}^{s}\{ -M_{l}-N_{l}, -M_{l}-N_{l}+1,\, \ldots, \, -M_{l} \}.$
In order to find the spectral characteristic $h_N(e^{i\lambda})$ of the estimate
\begin{equation*}
\hat{A}_N\vec{\xi}=\int\limits_{-\pi}^{\pi}(h_N(e^{i\lambda}))^\top(Z_{\xi}(d\lambda)+ Z_{\eta}(d\lambda))
 \end{equation*}
and the  mean-square error $\Delta(h_N;F,G)$ of the estimate  of the functional $A_N\vec{\xi}$, we define the vector $\vec{\bold{a}}_N^\top=(0, 0, \ldots,  0, \vec{a}_N^\top)$
which has first $T\cdot|S|=T\cdot\sum\limits_{k=1}^{s}(N_k+1)$ zero components and the last component is  $\vec{a}_N^\top=(\vec{a}(0)^\top,\vec{a}(1)^\top,\ldots, \vec{a}(N)^\top, 0, 0, \ldots).$

Consider the linear operator $\bold{R}_N$ in the space   $\ell_2$ which is defined as follows: $\bold{R}_N(k,j)=\bold{R}(k,j)$, $j\leq N$, $\bold{R}_N(k,j)=0$, $j> N$.

Thus the spectral characteristic of the optimal estimation $\hat{A}_N\xi$  can be calculated by the formula
\begin{equation}\label{4N} \begin{split}
(h_N(e^{i\lambda}))^\top=(A_N(e^{i\lambda}))^\top F(\lambda)(F(\lambda)+G(\lambda))^{-1}-\left(\sum\limits_{k\in U}(\bold{B}^{-1}\bold{R}_N\vec{\bold{a}}_N)(k) e^{ik\lambda}\right)^\top (F(\lambda)+G(\lambda))^{-1}.
\end{split} \end{equation}

The mean-square error of the estimate  $\hat{A}_N\vec{\xi}$  is defined by the formula
\[\Delta(h_N;F,G)=E\left|A_N\vec{\xi}-\hat{A}_N\vec{\xi}\right|^2=\frac{1}{2\pi}\int\limits_{-\pi}^{\pi}(r_G^N(\lambda))^\top F(\lambda)\overline{r_G^N(\lambda)}d\lambda+\frac{1}{2\pi}\int\limits_{-\pi}^{\pi}(r_F^N(\lambda))^\top G(\lambda)\overline{r_F^N(\lambda)}d\lambda=\]
\begin{equation} \label{55N}
=\langle\bold{R}_N\vec{\bold{a}}_N \bold{B}^{-1}\bold{R}_N\vec{\bold{a}}_N\rangle+\langle\bold{Q}_N\vec{\bold{a}}_N,\vec{\bold{a}}_N\rangle,
\end{equation}
where
\[(r_F^N(\lambda))^\top=\left((A_N(e^{i\lambda}))^\top F(\lambda)-\left(\sum\limits_{k \in U}(\bold{B}^{-1}\bold{R}_N\vec{\bold{a}_N})(k) e^{ik\lambda}\right)^\top\right)(F(\lambda)+G(\lambda))^{-1},\]
\[(r_G^N(\lambda))^\top=\left((A_N(e^{i\lambda}))^\top G(\lambda)+\left(\sum\limits_{k \in U}(\bold{B}^{-1}\bold{R}_N\vec{\bold{a}_N})(k) e^{ik\lambda}\right)^\top\right)(F(\lambda)+G(\lambda))^{-1},\]
and $\bold{Q}_N$ is the linear operator in the space  $\ell_2$, $\bold{Q}_N(k,j)=\bold{Q}(k,j)$, $k,j\leq N$, $\bold{Q}_N(k,j)=0$ if $k> N$ or $j> N$. Note, that linear operators $\bold{B}$, $\bold{R}$, $\bold{Q}$ are defined in Corollary \ref{c21}.

\begin{nas}\label{c23}
Let $\{\vec{\xi}(j), j\in \mathbb{Z}\}$ and $\{\vec{\eta}(j), j\in  \mathbb{Z}\}$ be uncorrelated multidimensional stationary stochastic sequences with spectral densities $F(\lambda)$ and $G(\lambda)$ which satisfy the minimality condition (\ref{minimal}). The spectral characteristic $h_N(e^{i\lambda})$ and the mean-square error  $\Delta(h_N;F,G)$ of the optimal linear estimate of the functional $A_N\vec{\xi}$ which depends on  unknown values of the sequence $\vec{\xi}(j)$ based on observations of the sequence $\vec{\xi}(j)+\vec{\eta}(j),$ $j\in \mathbb{Z}_{-}\backslash S$ can be calculated by formulas (\ref{4N}), (\ref{55N}).
\end{nas}

In the case where the sequence $\{\vec{\xi}(j), j\in \mathbb{Z}\}$ is observed without noise we have the following corollary.

\begin{nas} \label{c24}
Let $\{\vec{\xi}(j), j\in \mathbb{Z}\}$  be a multidimensional stationary stochastic sequence with the spectral density $F(\lambda)$ which satisfy the minimality condition (\ref{minimal1}). The spectral characteristic $h_N(e^{i\lambda})$ and the mean-square error  $\Delta(h_N,F)$ of the optimal linear estimate $\hat{A}_N\vec{\xi}$ of the functional $A_N\vec{\xi}$ can be calculated by the formulas   (\ref{spchar2N}), (\ref{6N})
 \begin {equation} \label{spchar2N}
(h_N(e^{i\lambda}))^\top=\left(\sum\limits_{j=0}^{N}\vec{a}(j)e^{ij\lambda }\right)^\top-\left(\sum\limits_{j \in U}\left(\bold{B}^{-1}\vec{\bold{a}}_N\right)(j) e^{ij\lambda }\right)^\top(F(\lambda))^{-1},
\end{equation}
\begin {equation} \label {6N} \begin{split}
 \Delta(h_N;F)=\langle\bold{B}^{-1}\vec{\bold{a}}_N, \vec{\bold{a}}_N\rangle.
\end{split} \end{equation}

The linear operator $\bold{B}$ is defined in Corollary \ref{c21}.
\end{nas}

In order to demonstrate the developed techniques we propose the following example.
%%%%%%%%%%%%%%%%%%%%%%%%%%%%%%%%%%%

\begin{exm} \label{example1}
Consider the problem of the optimal linear estimation of the functional
 $$A_1\vec{\xi} = \vec{a}(0)^\top\vec{\xi}(0) + \vec{a}(1)^\top\vec{\xi}(1)$$ which depends on the unknown values of a stationary sequence  $\vec{\xi}(j)=\{\xi_k(j)\}_{k=1}^2$ from  observations of the sequence $\vec{\xi}(j)$ at  points $j \in \mathbb{Z}_{-}\backslash S, $ where $S = \{-3, -2\}$, $\vec{a}(0)=(1, 1)^\top$, $\vec{a}(1)=(1, 1)^\top$. Let $\xi_{1} (n)= \xi (n)$ be a stationary stochastic sequence with the  spectral density $f( \lambda )$, and let $ \xi_{2} (n)= \xi (n)+ \eta (n)$, where $ \eta (n)$ is an uncorrelated with $ \xi (n)$ stationary stochastic sequence with the  spectral density $g( \lambda )$.  In this case the matrix of spectral densities is of the form
 \[F( \lambda )= \left( \begin{array}{cc} {f( \lambda )} & {f( \lambda )} \\{f( \lambda )} & {f( \lambda )+g( \lambda )} \end{array} \right), \]

\noindent and the inverse matrix is as follows
 \[(F( \lambda ))^{-1} = \left( \begin{array}{cc} { \frac{1}{f( \lambda )}+\frac{1}{g( \lambda )}} & { \frac{-1}{g( \lambda )}} \\{ \frac{-1}{g( \lambda )}} & { \frac{1}{g( \lambda )}} \end{array} \right). \]
Let
  \[f( \lambda )= \frac{1}{\left|1-b_{1} e^{i \lambda} \right|^{2}} , \quad g( \lambda )= \frac{1}{\left|1-b_{2} e^{i \lambda} \right|^{2}},\; b_{1},b_{2} \in R. \]
In this case the inverse matrix is of the form
\[(F( \lambda ))^{-1} = \left( \begin{array}{cc} { \left|1-b_1 e^{i \lambda} \right|^{2}+\left|1-b_2 e^{i \lambda} \right|^{2}} & {-\left|1-b_2 e^{i \lambda} \right|^{2}} \\{ -\left|1-b_2 e^{i \lambda} \right|^{2}} & { \left|1-b_2 e^{i \lambda} \right|^{2}} \end{array} \right)=B(-1)e^{-i\lambda} +B(0)+B(1)e^{i\lambda},\]
where \[B(0)=\left( \begin{array}{cc}{2 +b_1^2+b_2^2}&{-1-b_2^2}\\{-1-b_2^2}&{1+b_2^2}\end{array}\right),
\,B(1)=B(-1)=\left( \begin{array}{cc}{-b_1-b_2} & {b_2}\\{b_2}& {-b_2}\end{array}\right),\]
are the Fourier coefficients of the function $(F(\lambda))^{-1}$.

According to the Corollary \ref{c24} the spectral characteristic of the optimal estimate $\hat{A}_1\vec{\xi}$ of the functional $A_1\vec{\xi}$ is calculated by the formula
\begin{equation*}\begin{split}
(h_1(e^{i\lambda}))^\top &= \left(\vec{a}(0)+ \vec{a}(1)e^{i\lambda} \right)^\top - \left(\sum\limits_{j \in U}\left(\bold{B}^{-1}\vec{\bold{a}}_1\right)(j) e^{ij\lambda }\right)^\top\left(B(-1)e^{-i\lambda} +B(0)+B(1)e^{i\lambda}\right),
\end{split}\end{equation*}
where vector $\vec{\bold{a}}_1^\top = (0, 0, 0, 0, \vec{a}(0)^\top, \vec{a}(1)^\top, 0, 0,\ldots)$.

To find the unknown coefficients
 $\vec{c}(j)=\left(\bold{B}^{-1}\vec{\bold{a}}_1\right)(j), j \in U=S \cup \{0, 1, 2,\ldots\}$,
we use equation (\ref{rivn}), where $\vec{\bold{c}}^\top = (\vec{c}(-3)^\top,\vec{c}(-2)^\top,\vec{c}(0)^\top,\vec{c}(1)^\top,\vec{c}(2)^\top,\vec{c}(3)^\top, \ldots)$. The operator $\bold{B}$ is defined by the matrix
\begin{equation*}
B= \left(\begin{array}{cccccccc}
B(0) & B(-1) & 0 & 0 & 0 & 0& 0& \ldots \\
B(1) & B(0) & 0 & 0 & 0 & 0& 0& \ldots \\
0 & 0 & B(0) & B(-1) & 0 & 0& 0& \ldots \\
0 & 0 &  B(1) & B(0) & B(-1) & 0& 0& \ldots \\
0 & 0 & 0 & B(1) & B(0) &B(-1)& 0& \ldots \\
0 & 0 & 0 & 0 & B(1) & B(0)& B(-1)& \ldots \\
0 & 0 & 0 & 0 & 0 & B(1)& B(0)& \ldots \\
\ldots &   &   &   &  &   &   &
\end{array}\right).
\end{equation*}

We have to find the inverse matrix
 $B^{-1}$ which defines the inverse operator $\bold{B}^{-1}$. We first represent the  matrix $B$ in the form
\begin{equation*}
B= \left(\begin{array}{cc}
B_{00} & 0   \\
0 & B_{11}
\end{array}\right),
\end{equation*}
where
\begin{equation*}
B_{00}= \left(\begin{array}{ccccccc}
  B(0) & B(-1)\\
  B(1) & B(0)
  \end{array}\right),
\end{equation*}
\begin{equation*}
B_{11}= \left(\begin{array}{cccccc}
  B(0) & B(-1) & 0 & 0& 0& \ldots \\
   B(1) & B(0) & B(-1) & 0& 0& \ldots \\
0 & B(1) & B(0) & B(-1) & 0& \ldots \\
  0 & 0 & B(1) & B(0) & B(-1) & \ldots \\
   0 & 0 & 0 & B(1) & B(0) & \ldots \\
\ldots   &   &   &  &   &
\end{array}\right).
\end{equation*}
Making use of the indicated representation we may conclude that the  matrix $B^{-1}$ can be represented in the form
\begin{equation*}
B^{-1}= \left(\begin{array}{cc}
B_{0 0}^{-1} & 0   \\
0 & B_{1 1}^{-1}
\end{array}\right),
\end{equation*}
where $B_{0 0}^{-1}$, $B_{1 1}^{-1}$ are inverse matrices to the matrices  $B_{0 0} $, $B_{1 1} $ respectively.
The matrix $B_{0 0}^{-1}$ can be found in the form
 \[B_{0 0}^{-1}=\left( \begin{array}{cccccc} {\frac{1+b_1^2}{A}} & {\frac{1+b_1^2}{A}} & {\frac{b_1}{A}} & {\frac{b_1}{A}}\\ \frac{1+b_1^2}{A}& \frac{1+b_1^2}{A}+\frac{1+b_2^2}{B} & {\frac{b_1}{A}} & {\frac{b_1}{A}}+{\frac{b_2}{B}}\\ {\frac{b_1}{A}} & {\frac{b_1}{A}} & {\frac{1+b_1^2}{A}} & {\frac{1+b_1^2}{A}}\\{\frac{b_1}{A}} & {\frac{b_1}{A}}+{\frac{b_2}{B}} & {\frac{1+b_1^2}{A}} & {\frac{1+b_1^2}{A}}+{\frac{1+b_2^2}{B}} \end{array} \right),\]
where $A=1+b_1^2+b_1^4$, $B=1+b_2^2+b_2^4$.
In order to find the matrix $(B_{1 1})^{-1}$ we use the following method. The matrix $B_{1 1}$ is constructed with the help of the Fourier coefficients of the function   $(F(\lambda))^{-1}$
\begin{equation*}\label{matr22}
  B_{1 1}(k,j)= B(k-j), \quad k,j= 0, 1, 2,\ldots.
 \end{equation*}
The density $(F(\lambda))^{-1}$ admits the factorization
\[
(F(\lambda))^{-1}=\sum\limits_{p=-\infty}^{\infty}B(p)e^{ip\lambda}=\left(\sum\limits_{j=0}^{\infty}\psi(j)e^{-ij\lambda}\right)\cdot\left(\sum\limits_{j=0}^{\infty}\psi(j)e^{-ij\lambda}\right)^*=\]
\[=\left(\sum\limits_{j=0}^{\infty}\theta(j)e^{-ij\lambda}\right)^*\cdot\left(\sum\limits_{j=0}^{\infty}\theta(j)e^{-ij\lambda}\right),\]
where
\[\psi(0)=\left( \begin{array}{cc} {1} & {1} \\{0} & {-1} \end{array} \right), \psi(1) =\left( \begin{array}{cc} {-b_1} & {-b_2} \\{0} & {b_2} \end{array} \right),  \psi(j)=0, j>1,\; \theta(j)= \left( \begin{array}{cc} {b_1^j} & {b_1^j} \\{0} & {-b_2^j} \end{array} \right),  j\geq 0.\]
Hence $B(p)=\sum\limits_{k=0}^{\infty}\psi(k)(\psi(k+p))^*,$ $p\geq0$, and $B(-p)=(B(p))^*$, $p\geq0$. Then $$B(i-j) = \sum\limits_{l=\max(i,j)}^{\infty}\psi(l-i)(\psi(l-j))^*.$$
Denote by $\bold{\Psi}$ and $ \bold{\Theta}$  linear operators in the space  $\ell_2$ determined by  matrices with elements  $\bold{\Psi}(i,j)=\psi(j-i)$, $\bold{\Theta}(i,j)=\theta(j-i)$, for $0\leq i\leq j$,  $\bold{\Psi}(i,j)=0$, $\bold{\Theta}(i,j)=0$, for  $0\leq j< i$. Then elements of the matrix $B_{11}$ can be represented in the form
$B_{1 1}(i,j)=(\bold{\Psi}\bold{\Psi}^*)(i,j).$
Since  $\bold{\Psi}\bold{\Theta}=\bold{\Theta}\bold{\Psi}=I$, elements of the matrix $B_{1 1}^{-1}$ can be calculated by the formula $B_{1 1}^{-1}(i,j)=(\bold{\Theta}^*\bold{\Theta})(i,j)=\sum\limits_{l=0}^{\min(i,j)}(\theta(i-l))^*\theta(j-l)$.

From  equation (\ref{rivn}) we can find the unknown coefficients  $c(j), j \in U$,
\begin{equation*}\begin{split}
&\vec{c}(-3)= \vec{0},\\
&\vec{c}(-2)= \vec{0},\\
&\vec{c}(0)= B_{11}^{-1}(0,0)\vec{a}(0)+B_{11}^{-1}(0,1)\vec{a}(1),\\
&\vec{c}(1)=B_{11}^{-1}(1,0)\vec{a}(0)+B_{11}^{-1}(1,1)\vec{a}(1),\\
&\vec{c}(2)=B_{11}^{-1}(2,0)\vec{a}(0)+B_{11}^{-1}(2,1)\vec{a}(1),\\
&\ldots\\
&\vec{c}(i)=B_{11}^{-1}(i,0)\vec{a}(0)+B_{11}^{-1}(i,1)\vec{a}(1), \quad i>2.
\end{split}\end{equation*}
Hence the spectral characteristic of the optimal estimate is calculated by the formula
\[
(h_1(e^{i\lambda}))^\top= \left(\vec{a}(0) + \vec{a}(1)e^{i\lambda} \right)^\top -
(\vec{c}(-3)e^{-i3\lambda}+\vec{c}(-2)e^{-i2\lambda}+\vec{c}(0)+\vec{c}(1)e^{i\lambda}+\vec{c}(2)e^{i2\lambda}+\]
\[+\sum\limits_{j > 2}\vec{c}(j) e^{ij\lambda }) \left(B(-1)e^{-i\lambda} +B(0)+B(1)e^{i\lambda}\right)
=-\vec{c}(0)^\top B(-1)e^{-i\lambda}-\vec{c}(1)^\top B(1)e^{ i2\lambda}-\]
\[-\vec{c}(2)^\top B(0)e^{i2\lambda}-\vec{c}(2)^\top B(1)e^{i3\lambda}-\sum\limits_{j > 2}\vec{c}(j)^\top e^{ij\lambda }\left(B(-1) e^{-i\lambda} +B(0)+B(1) e^{i\lambda}\right).
\]
Since coefficients $\vec{c}(j-1)^\top B(1)+\vec{c}(j)^\top B(0)+\vec{c}(j+1)^\top B(-1)$  for $j\geq 2$ are  zero, the spectral characteristic of the estimate  $\hat{A}_1\vec{\xi}$ is of the form
\[
(h_1(e^{i\lambda})^\top=-\vec{c}(0)^\top B(-1)e^{-i\lambda}=\left(b_2+b_2^2-2(b_1+b_1^2), \; -b_2-b_2^2\right)e^{-i\lambda}.
\]
The mean-square error of the estimate of the functional $A_1\vec{\xi}$ is calculated by the formula
\begin{equation*} \begin{split}
 \Delta(h_1;F)&=\langle B^{-1}\vec{\bold{a}}_1, \vec{\bold{a}}_1 \rangle=10+8b_1+4b_1^2+2b_2+b_2^2.
\end{split} \end{equation*}
\end{exm}

%------------------------------------------------------------------------------------------------

\section{Minimax approach to extrapolation problem for stationary sequences with missing observations}

Theorem \ref{t2}  and its Corollaries \ref{c21} -- \ref{c24} can be applied for finding solutions of the extrapolation problem for multidimensional stationary sequences with missing observations only in the  case  of spectral certainty, where spectral densities $F(\lambda),F_{\xi\eta}(\lambda),F_{\eta\xi}(\lambda),G(\lambda)$ are exactly known. If the complete information about spectral densities is impossible while a class of admissible spectral density matrices $D$ is given, the minimax(robust) method of extrapolation is reasonable. It consists in finding an estimate which minimizes the value of the mean-square error for all spectral density matrices from the given class of densities.  For description of the minimax method we introduce the following definitions (see
Moklyachuk \cite{Moklyachuk:2008,Moklyachuk:2015}, and Moklyachuk and Masytka \cite{Moklyachuk:Mas2008} - \cite{Moklyachuk:Mas2012}).

 \begin{ozn}
  For a given class of spectral density matrices $D$ the spectral densities $(F^0(\lambda),F^0_{\xi\eta}(\lambda),F^0_{\eta\xi}(\lambda),G^0(\lambda))\in D$ are called the least favorable in the class  $D$ for the optimal linear extrapolation of the functional $A\vec{\xi}$ if the following relation holds true
   $$\Delta\left(h\left(F^0,F^0_{\xi\eta},F^0_{\eta\xi},G^0\right);F^0,F^0_{\xi\eta},F^0_{\eta\xi},G^0\right)=\max\limits_{(F,F_{\xi\eta},F_{\eta\xi},G)\in D}\Delta\left(h\left(F,F_{\xi\eta},F_{\eta\xi},G\right);F,F_{\xi\eta},F_{\eta\xi},G\right).$$
\end{ozn}

\begin{ozn}
For a given class of spectral density matrices $D$  the spectral characteristic  $h^0(e^{i\lambda})$  of the optimal linear estimate  of the functional $A\vec{\xi}$ is called minimax-robust if there are satisfied conditions
$$h^0(e^{i\lambda})\in H_{ D}= \bigcap\limits_{(F,F_{\xi\eta},F_{\eta\xi},G)\in D} L_2^s(F+G),$$
$$\min\limits_{h\in H_{ D}}\max\limits_{{(F,F_{\xi\eta},F_{\eta\xi},G)\in  D}}\Delta\left(h;F,F_{\xi\eta},F_{\eta\xi},G\right)=\max\limits_{{(F,F_{\xi\eta},F_{\eta\xi},G)\in  D}}\Delta\left(h^0;F,F_{\xi\eta},F_{\eta\xi},G\right).$$
\end{ozn}

From the introduced definitions and formulas derived above we can obtain the following statement.
%\begin{comment}
\begin{lem}
Spectral densities $(F^0(\lambda),F^0_{\xi\eta}(\lambda),F^0_{\eta\xi}(\lambda),G^0(\lambda))\in D$, satisfying the minimality condition (\ref{minimal}), are the least favorable in the class $ D$ for the optimal linear extrapolation of the functional $A\vec{\xi}$ if the Fourier coefficients \eqref{Fcoef} of the functions $$(F^0_{\zeta}(\lambda))^{-1}, \quad (F^0(\lambda)+ F_{\xi\eta}^0(\lambda))(F^0_{\zeta}(\lambda))^{-1}, \quad F^0(\lambda)(F^0_{\zeta}(\lambda))^{-1}G^0(\lambda)-F_{\xi\eta}^0(\lambda)(F^0_{\zeta}(\lambda))^{-1}F_{\eta\xi}^0(\lambda)$$ define operators $\bold{B}^0, \bold{R}^0, \bold{Q}^0$ which determine a solution of the constrained optimization problem
\begin{equation} \label{extrem}
\max\limits_{{(F,F_{\xi\eta},F_{\eta\xi},G)\in D}}\,(\langle\bold{R}\vec{\bold{a}},\bold{B}^{-1}\bold{R}\vec{\bold{a}}\rangle +\langle\bold{Q}\vec{\bold{a}},\vec{\bold{a}}\rangle)=\langle\bold{R}^0\vec{\bold{a}},(\bold{B}^0)^{-1}\bold{R}^0\vec{\bold{a}}\rangle+\langle\bold{Q}^0\vec{\bold{a}},\vec{\bold{a}}\rangle.
\end{equation}
The minimax spectral characteristic $h^0(e^{i\lambda})=h(F^0,F^0_{\xi\eta},F^0_{\eta\xi},G^0)$ is calculated by the formula (\ref{4}) if $h(F^0,F^0_{\xi\eta},F^0_{\eta\xi},G^0) \in H_{D}.$
\end{lem}

In the case of uncorrelated stationary sequences  the corresponding definitions and lemmas are as follows.

 \begin{ozn}
  For a given class of spectral densities $D=D_F \times D_G$ the spectral densities $F^0(\lambda) \in D_F$, $G^0(\lambda) \in D_G$ are called the least favorable in the class  $D$ for the optimal linear extrapolation of the functional $A\vec{\xi}$ based on observations of the uncorrelated sequences if the following relation holds true
   $$\Delta\left(F^0,G^0\right)=\Delta\left(h\left(F^0,G^0\right);F^0,G^0\right)=\max\limits_{(F,G)\in D_F\times D_G}\Delta\left(h\left(F,G\right);F,G\right).$$
\end{ozn}

\begin{ozn}
For a given class of spectral densities $D=D_F\times D_G$ the spectral characteristic  $h^0(e^{i\lambda})$  of the optimal linear estimate  of the functional $A\vec{\xi}$ based on observations of the uncorrelated sequences is called minimax-robust if there are satisfied conditions
$$h^0(e^{i\lambda})\in H_D= \bigcap\limits_{(F,G)\in D_F\times D_G} L_2^s(F+G),$$
$$\min\limits_{h\in H_D}\max\limits_{(F,G)\in D}\Delta\left(h;F,G\right)=\max\limits_{(F,G)\in D}\Delta\left(h^0;F,G\right).$$
\end{ozn}

\begin{lem}
Spectral densities $F^0(\lambda)\in D_F,$ $G^0(\lambda) \in D_G$ satisfying the minimality condition (\ref{minimal-Unc}) are the least favorable in the class $D=D_F\times D_G$ for the optimal linear extrapolation of the functional $A\vec{\xi}$ based on observations of the uncorrelated sequences if the Fourier coefficients \eqref{Fcoef-Un} of functions $$(F^0(\lambda)+G^0(\lambda))^{-1}, \quad F^0(\lambda)(F^0(\lambda)+G^0(\lambda))^{-1}, \quad F^0(\lambda)(F^0(\lambda)+G^0(\lambda))^{-1}G^0(\lambda)$$ define operators $\bold{B}^0, \bold{R}^0, \bold{Q}^0$ which determine a solution of the constrained optimization problem
\begin{equation} \label{extrem}
\max\limits_{(F,G)\in D_F\times D_G}(\langle\bold{R}\vec{\bold{a}},\bold{B}^{-1}\bold{R}\vec{\bold{a}}\rangle +\langle\bold{Q}\vec{\bold{a}},\vec{\bold{a}}\rangle)= \langle\bold{R}^0\vec{\bold{a}},(\bold{B}^0)^{-1}\bold{R}^0\vec{\bold{a}}\rangle+\langle\bold{Q}^0\vec{\bold{a}},\vec{\bold{a}}\rangle.
\end{equation}
The minimax spectral characteristic $h^0=h(F^0,G^0)$ is calculated by the formula (\ref{sp-ch2}) if $h(F^0,G^0) \in H_D.$
\end{lem}

In the case of observations of the sequence without noise we obtain the following corollary.

  \begin{nas}
Let the spectral density  $F^0(\lambda)\in D_F$  be such that the function $(F^0(\lambda))^{-1}$ satisfies the minimality condition. The spectal density $F^0(\lambda)\in D_F$ is the least favorable in the class   $ D_F$ for the optimal linear extrapolation of the functional $A\vec{\xi}$ if the Fourier coefficients of the function $(F^0(\lambda))^{-1}$ define the operator $\bold{B}^0$ which determines a solution of the optimization problem
 \begin{equation} \label{extrem2}
\max\limits_{F\in D_F}\langle\bold{B}^{-1}\vec{\bold{a}},\vec{\bold{a}}\rangle=\langle(\bold{B}^0)^{-1}\vec{\bold{a}},\vec{\bold{a}}\rangle.
\end{equation}
The minimax spectral characteristic $h^0=h(F^0)$ is calculated by the formula (\ref{spchar2}) if $h(F^0) \in H_{D_F}.$
 \end{nas}

The least favorable spectral densities $F^0(\lambda)$, $G^0(\lambda)$ and the minimax spectral characteristic $h^0=h(F^0,G^0)$ form a saddle point of the function $\Delta \left(h;F,G\right)$ on the set $H_D\times D.$ The saddle point inequalities
$$\Delta\left(h;F^0,G^0\right)\geq\Delta\left(h^0;F^0,G^0\right)\geq \Delta\left(h^0;F,G\right) $$  $$ \forall h \in H_D, \forall F \in D_F, \forall G \in D_G$$
hold true if $h^0=h(F^0,G^0)$ and $h(F^0,G^0)\in H_D,$  where $(F^0,G^0)$ is a solution of the constrained optimization problem
\begin{equation} \label{7}
\sup\limits_{(F,G)\in D_F\times D_G}\Delta\left(h(F^0,G^0);F,G\right)=\Delta\left(h(F^0,G^0);F^0,G^0\right),
\end{equation}

\begin{equation*}
\Delta\left(h(F^0,G^0);F,G\right)=\frac{1}{2\pi}\int\limits_{-\pi}^{\pi}(r_G^0(\lambda))^\top F(\lambda)\overline{r_G^0(\lambda)}d\lambda+\frac{1}{2\pi}\int\limits_{-\pi}^{\pi}(r_F^0(\lambda))^\top G(\lambda)\overline{r_F^0(\lambda)}d\lambda,
\end{equation*}
\[(r_F^0(\lambda))^\top=\left((A(e^{i\lambda}))^\top F^0(\lambda)-\left(\sum\limits_{k \in U}((\bold{B}^0)^{-1}\bold{R}^0\vec{\bold{a}})(k) e^{ik\lambda}\right)^\top\right)(F^0(\lambda)+G^0(\lambda))^{-1},\]
\[(r_G^0(\lambda))^\top=\left((A(e^{i\lambda}))^\top G^0(\lambda)+\left(\sum\limits_{k \in U}((\bold{B}^0)^{-1}\bold{R}^0\vec{\bold{a}})(k) e^{ik\lambda}\right)^\top\right)(F^0(\lambda)+G^0(\lambda))^{-1}.\]

The constrained optimization problem  (\ref{7}) is equivalent to the unconstrained optimization problem (see Pshenichnyj \cite{Pshenychn}):
\begin{equation} \label{8}
\Delta_D(F,G)=-\Delta(h(F^0,G^0);F,G)+\delta((F,G)\left|D_F\times D_G\right.)\rightarrow \inf,
\end{equation}
 where $\delta((F,G)\left|D_F\times D_G\right.)$ is the indicator function of the set $D=D_F\times D_G$.

A solution of the problem   (\ref{8}) is determined by the condition $0 \in \partial\Delta_D(F^0,G^0),$
where $\partial\Delta_D(F^0,G^0)$ is the subdifferential of the convex functional $\Delta_D(F,G)$ at point $(F^0,G^0)$. This condition
is the necessary and sufficient condition under which the pair $(F^0,G^0)$ belongs to the set of minimums of the convex functional $\Delta(h(F^0,G^0);F,G)$.
This condition makes it possible to find the least favourable spectral densities in some special classes of spectral densities $D$ (see books by Ioffe and Tihomirov \cite{Ioffe}, Pshenichnyj \cite{Pshenychn}, Rockafellar \cite{Rockafellar}).

Note, that the form of the functional $\Delta\left(h^0;F,G\right)$ is convenient for application the Lagrange method of indefinite multipliers for finding solution to the problem (\ref{8}).
Making use the method of Lagrange multipliers and the form of
subdifferentials of the indicator functions we describe relations that determine the least favourable spectral densities in some special classes of spectral densities (see books by Moklyachuk \cite{Moklyachuk:2008b,Moklyachuk:2008}, Moklyachuk and Masyutka \cite{Moklyachuk:Mas2012} for additional details).

\section{Least favorable spectral densities in the class $D=D_0 \times D_V^U$}
Consider the problem of extrapolation of the functional $A\vec{\xi}$  based on observations of the uncorrelated sequences in the case where spectral densities of the observed sequences belong to the class  $D=D_0 \times D_V^U$, where
$$ D_{0}^{1} =\bigg\{F(\lambda )\left|\frac{1}{2\pi }
\int _{-\pi }^{\pi}{\rm{Tr}}\, F(\lambda )d\lambda =p\right.\bigg\},$$
$${D_{V}^{U}} ^{1}  =\bigg\{G(\lambda )\bigg|{\mathrm{Tr}}\, V(\lambda
)\le {\mathrm{Tr}}\, G(\lambda )\le {\mathrm{Tr}}\, U(\lambda ), \frac{1}{2\pi } \int _{-\pi}^{\pi}{\mathrm{Tr}}\,  G(\lambda)d\lambda  =q \bigg\},$$
$$D_{0}^{2} =\bigg\{F(\lambda )\left|\frac{1}{2\pi }
\int _{-\pi}^{\pi}f_{kk} (\lambda )d\lambda =p_{k}, k=\overline{1,T}\right.\bigg\},$$
$${D_{V}^{U}} ^{2}  =\bigg\{G(\lambda )\bigg|v_{kk} (\lambda )  \le
g_{kk} (\lambda )\le u_{kk} (\lambda ), \frac{1}{2\pi} \int _{-\pi}^{\pi}g_{kk} (\lambda
)d\lambda  =q_{k} , k=\overline{1,T}\bigg\},$$
$$D_{0}^{3} =\bigg\{F(\lambda )\left|\frac{1}{2\pi} \int _{-\pi}^{\pi}\left\langle B_{1} ,F(\lambda )\right\rangle d\lambda  =p\right.\bigg\},$$
$${D_{V}^{U}} ^{3}  =\bigg\{G(\lambda )\bigg|\left\langle B_{2}
,V(\lambda )\right\rangle \le \left\langle B_{2} ,G(\lambda
)\right\rangle \le \left\langle B_{2} ,U(\lambda)\right\rangle,\frac{1}{2\pi }
\int _{-\pi}^{\pi}\left\langle B_{2},G(\lambda)\right\rangle d\lambda  =q\bigg\},$$
$$D_{0}^{4} =\bigg\{F(\lambda )\left|\frac{1}{2\pi} \int
_{-\pi}^{\pi}F(\lambda )d\lambda  =P\right.\bigg\},$$
$${D_{V}^{U}} ^{4}=\left\{G(\lambda )\bigg|V(\lambda )\le G(\lambda
)\le U(\lambda ), \frac{1}{2\pi } \int _{-\pi}^{\pi}G(\lambda )d\lambda=Q\right\}.$$
Here spectral densities $V( \lambda ),U( \lambda )$ are known and fixed, $p, q, p_k, q_k, k=\overline{1,T}$ are given numbers, $P, Q, B_1, B_2$ are given positive definite Hermitian matrices.

From the condition $0\in \partial \Delta _{D} (F^{0} ,G^{0} )$ we find the following equations which determine the least favourable spectral densities for these given sets of admissible spectral densities.

For the first pair $D_{0}^{1}\times {D_{V}^{U}} ^{1}$ we have equations
\begin{equation} \label{eq_4_1}
(r_G^0(\lambda))^{*}(r_G^0(\lambda))^\top=\alpha^{2} (F^{0} (\lambda )+G^{0} (\lambda
))^{2} ,
\end{equation}
\begin{equation} \label{eq_4_2}
(r_F^0(\lambda))^{*}(r_F^0(\lambda))^\top=(\beta^{2} +\gamma _{1}
(\lambda )+\gamma _{2} (\lambda ))(F^{0} (\lambda )+G^{0} (\lambda))^{2},
\end{equation}
where $\alpha^{2}, \beta^{2}$ are Lagrange multipliers,  $\gamma _{1} (\lambda )\le 0$ and $\gamma _{1} (\lambda )=0$ if ${\mathrm{Tr}}\,
G^{0} (\lambda )> {\mathrm{Tr}}\,  V(\lambda ),$ $\gamma _{2} (\lambda )\ge 0$ and $\gamma _{2} (\lambda )=0$ if $ {\mathrm{Tr}}\,G^{0}(\lambda )< {\mathrm{Tr}}\,  U(\lambda).$

For the second pair $D_{0}^{2}\times {D_{V}^{U}} ^{2}$ we have equations
\begin{equation}  \label{eq_4_3}
(r_G^0(\lambda))^{*}(r_G^0(\lambda))^\top=(F^{0} (\lambda )+G^{0} (\lambda
))\left\{\alpha _{k}^{2} \delta _{kl} \right\}_{k,l=1}^{T} (F^{0}
(\lambda )+G^{0} (\lambda )),
\end{equation}
\begin{equation}\label{eq_4_4}
(r_F^0(\lambda))^{*}(r_F^0(\lambda))^\top=(F^{0} (\lambda)+G^{0} (\lambda ))\left\{(\beta_{k}^{2} +\gamma _{1k} (\lambda )+\gamma _{2k}(\lambda ))\delta _{kl} \right\}_{k,l=1}^{T} (F^{0} (\lambda )+G^{0}(\lambda )),
\end{equation}
where  $\alpha _{k}^{2}, \beta_{k}^{2}$ are Lagrange multipliers, $\delta _{kl}$ are Kronecker symbols, $\gamma _{1k} (\lambda )\le 0$ and $\gamma _{1k} (\lambda )=0$ if $g_{kk}^{0} (\lambda )>v_{kk} (\lambda ),$ $\gamma _{2k} (\lambda )\ge 0$ and $\gamma _{2k} (\lambda )=0$ if $g_{kk}^{0} (\lambda )<u_{kk} (\lambda).$

For the third pair $D_{0}^{3}\times {D_{V}^{U}} ^{3}$ we have equations
\begin{equation}  \label{eq_4_5}
(r_G^0(\lambda))^{*}(r_G^0(\lambda))^\top=\alpha^{2} (F^{0} (\lambda )+G^{0} (\lambda
))B_{1}^{\top} (F^{0} (\lambda )+G^{0} (\lambda )),
\end{equation}
\begin{equation}\label{eq_4_6}
(r_F^0(\lambda))^{*}(r_F^0(\lambda))^\top=(\beta^{2} +\gamma'_{1} (\lambda )+\gamma'_{2} (\lambda
))(F^{0} (\lambda )+G^{0} (\lambda ))B_{2}^{\top}(F^{0} (\lambda)+G^{0} (\lambda )),
\end{equation}
where $\alpha^{2}, \beta^{2}$ are Lagrange multipliers, $\gamma'_{1}( \lambda )\le 0$ and $\gamma'_{1} ( \lambda )=0$ if $\langle B_{2},G^{0} ( \lambda \rangle > \langle B_{2},V( \lambda ) \rangle,$ $\gamma'_{2}( \lambda )\ge 0$ and $\gamma'_{2} ( \lambda )=0$ if $\langle
B_{2} ,G^{0} ( \lambda \rangle < \langle B_{2} ,U( \lambda ) \rangle.$

For the fourth pair $D_{0}^{4}\times {D_{V}^{U}} ^{4}$ we have equations
\begin{equation} \label{eq_4_7}
(r_G^0(\lambda))^{*}(r_G^0(\lambda))^\top=(F^{0} (\lambda )+G^{0} (\lambda
))\vec{\alpha}\cdot \vec{\alpha}^{*}(F^{0} (\lambda )+G^{0}
(\lambda )),
\end{equation}
\begin{equation}\label{eq_4_8}
(r_F^0(\lambda))^{*}(r_F^0(\lambda))^\top=(F^{0} (\lambda )+G^{0} (\lambda ))(\vec{\beta}\cdot \vec{\beta}^{*}+\Gamma _{1} (\lambda )+\Gamma _{2} (\lambda ))(F^{0} (\lambda)+G^{0} (\lambda ))
\end{equation}
where $\vec{\alpha}, \vec{\beta}$ are Lagrange multipliers, $\Gamma _{1} (\lambda )\le 0$ and $\Gamma _{1} (\lambda )=0$ if $G^{0}(\lambda )>V(\lambda ),$ $
\Gamma _{2} (\lambda )\ge 0$ and $\Gamma _{2} (\lambda )=0$ if $G^{0}(\lambda )<U(\lambda ).$

The following theorem and corollaries hold true.

\begin{thm}
Let the minimality condition (\ref{minimal-Unc}) hold true. The least favorable spectral densities  $F^0(\lambda)$, $G^0(\lambda)$  in the classes $D_0 \times D_V^U$ for the optimal linear extrapolation of the functional $A\vec{\xi}$ are determined by relations
(\ref{eq_4_1}), (\ref{eq_4_2}) for the first pair $D_{0}^{1}\times {D_{V}^{U}} ^{1}$ of sets of admissible spectral densities;
(\ref{eq_4_3}), (\ref{eq_4_4}) for the second pair $D_{0}^{2}\times {D_{V}^{U}} ^{2}$ of sets of admissible spectral densities;
(\ref{eq_4_5}), (\ref{eq_4_6}) for the third pair $D_{0}^{3}\times {D_{V}^{U}} ^{3}$ of sets of admissible spectral densities;
(\ref{eq_4_7}), (\ref{eq_4_8}) for the fourth pair $D_{0}^{4}\times {D_{V}^{U}} ^{4}$ of sets of admissible spectral densities;
constrained optimization problem (\ref{extrem}) and restrictions  on densities from the corresponding classes $D_0 \times D_V^U$.  The minimax-robust spectral characteristic of the optimal estimate of the functional $A\vec{\xi}$ is determined by the formula (\ref{sp-ch2}).
\end{thm}

\begin{nas}
Let the minimality condition (\ref{minimal1}) hold true. The least favorable spectral densities $F^{0}(\lambda)$ in the classes $D_0^{k}$, $k=1,2,3,4$, for the optimal linear extrapolation of the functional  $A\vec{\xi}$ from observations of the sequence $\vec{\xi}(j)$ at points  $j\in\mathbb{Z}_{-}\backslash S $, where $S=\bigcup\limits_{l=1}^{s}\{-M_{l}-N_l, \ldots,-M_{l} \},$ are determined by the following  equations, respectively,
\begin{equation}
((C^{0}(\lambda) )^{\top} )^{*}\cdot(C^{0}(\lambda) )^{\top}=\alpha^{2}(F^{0} (\lambda ))^{2},
\end{equation}
\begin{equation}
((C^{0}(\lambda) )^{\top} )^{*}\cdot(C^{0}(\lambda) )^{\top}=F^{0} (\lambda )\left\{\alpha _{k}^{2}\delta _{kl} \right\}_{k,l=1}^{T}F^{0} (\lambda ),
\end{equation}
\begin{equation}
((C^{0}(\lambda) )^{\top} )^{*}\cdot(C^{0}(\lambda) )^{\top}=\alpha^{2}F^{0} (\lambda )B_1^\top F^{0} (\lambda ),
\end{equation}
\begin{equation}
((C^{0}(\lambda) )^{\top} )^{*}\cdot(C^{0}(\lambda) )^{\top}=F^{0} (\lambda )\vec{\alpha}\cdot \vec{\alpha}^{*}F^{0} (\lambda ),
\end{equation}
constrained optimization problem (\ref{extrem2}) and restrictions  on densities from the corresponding classes  $D_0^{k}$, $k=1,2,3,4$. The minimax spectral characteristic of the optimal estimate of the functional $A\vec{\xi}$ is determined by the formula (\ref{spchar2}).
\end{nas}

\begin{nas}
Let the minimality condition (\ref{minimal1}) hold true. The least favorable spectral densities $F^{0}(\lambda)$ in the classes ${D_{V}^{U}} ^{k}$, $k=1,2,3,4$, for the optimal linear extrapolation of the functional  $A\vec{\xi}$ from observations of the sequence $\vec{\xi}(j)$ at points  $j\in\mathbb{Z}_{-}\backslash S $, where $S=\bigcup\limits_{l=1}^{s}\{-M_{l}-N_l, \ldots,-M_{l} \},$ are determined by the following  equations, respectively,
\begin{equation}
((C^{0}(\lambda) )^{\top} )^{*}\cdot(C^{0}(\lambda) )^{\top}=(\beta^{2} +\gamma _{1} (\lambda )+\gamma _{2} (\lambda )) (F^{0} (\lambda ))^{2},
\end{equation}
\begin{equation}
((C^{0}(\lambda) )^{\top} )^{*}\cdot(C^{0}(\lambda) )^{\top}=F^{0} (\lambda )\left\{(\beta_{k}^{2} +\gamma _{1k} (\lambda )+\gamma _{2k} (\lambda ))\delta _{kl}\right\}_{k,l=1}^{T}F^{0} (\lambda ),
\end{equation}
\begin{equation}
((C^{0}(\lambda) )^{\top} )^{*}\cdot(C^{0}(\lambda) )^{\top}=(\beta^{2} +\gamma'_{1}(\lambda )+\gamma'_{2}(\lambda )) F^{0} (\lambda )B_2^\top F^{0} (\lambda ),
\end{equation}
\begin{equation}
((C^{0}(\lambda) )^{\top} )^{*}\cdot(C^{0}(\lambda) )^{\top}=F^{0} (\lambda )(\vec{\beta}\cdot \vec{\beta}^{*}+\Gamma _{1} (\lambda )+\Gamma _{2} (\lambda ))F^{0} (\lambda ),
\end{equation}
constrained optimization problem (\ref{extrem2}) and restrictions  on densities from the corresponding classes  ${D_{V}^{U}} ^{k}$, $k=1,2,3,4$. The minimax spectral characteristic of the optimal estimate of the functional $A_s\vec{\xi}$ is determined by the formula (\ref{spchar2}).
\end{nas}

\section{Least favorable spectral densities in the class $D=D_{\varepsilon}\times D_{1\delta}$}
Consider the problem of minimax  extrapolation of the functional $A\vec{\xi}$  based on observations of the uncorrelated sequences in the case where spectral densities of the observed sequences belong to the class
 $D=D_{\varepsilon}\times D_{1\delta}$,
$$D_{\varepsilon }^{1}  =\bigg\{F(\lambda )\bigg|{\mathrm{Tr}}\,
F(\lambda )=(1-\varepsilon ) {\mathrm{Tr}}\,  F_{1} (\lambda
)+\varepsilon {\mathrm{Tr}}\,  W(\lambda ), \frac{1}{2\pi} \int _{-\pi}^{\pi}{\mathrm{Tr}}\,
F(\lambda )d\lambda =p \bigg\};$$
$$D_{1\delta}^{1}=\left\{G(\lambda )\biggl|\frac{1}{2\pi} \int_{-\pi}^{\pi}\left|{\rm{Tr}}(G(\lambda )-G_{1} (\lambda))\right|d\lambda \le \delta\right\};$$
$$D_{\varepsilon }^{2}  =\bigg\{F(\lambda )\bigg|f_{kk} (\lambda)
=(1-\varepsilon )f_{kk}^{1} (\lambda )+\varepsilon w_{kk}(\lambda), \frac{1}{2\pi} \int _{-\pi}^{\pi}f_{kk} (\lambda)d\lambda  =p_{k} , k=\overline{1,T}\bigg\};$$
$$D_{1\delta}^{2}=\left\{G(\lambda )\biggl|\frac{1}{2\pi } \int_{-\pi}^{\pi}\left|g_{kk} (\lambda )-g_{kk}^{1} (\lambda)\right|d\lambda  \le \delta_{k}, k=\overline{1,T}\right\};$$
$$D_{\varepsilon }^{3} =\bigg\{F(\lambda )\bigg|\left\langle B_{1},F(\lambda )\right\rangle =(1-\varepsilon )\left\langle B_{1},F_{1} (\lambda )\right\rangle+\varepsilon \left\langle B_{1},W(\lambda )\right\rangle, \frac{1}{2\pi}\int _{-\pi}^{\pi}\left\langle B_{1} ,F(\lambda )\right\rangle d\lambda =p\bigg\};$$
$$D_{1\delta}^{3}=\left\{G(\lambda )\biggl|\frac{1}{2\pi } \int_{-\pi}^{\pi}\left|\left\langle B_{2} ,G(\lambda )-G_{1}(\lambda )\right\rangle \right|d\lambda  \le \delta\right\};$$
$$D_{\varepsilon }^{4}=\left\{F(\lambda )\bigg|F(\lambda)=(1-\varepsilon )F_{1} (\lambda )+\varepsilon W(\lambda ),
\frac{1}{2\pi } \int _{-\pi}^{\pi}F(\lambda )d\lambda=P\right\},$$
$$D_{1\delta}^{4}=\left\{G(\lambda )\biggl|\frac{1}{2\pi} \int_{-\pi}^{\pi}\left|g_{ij} (\lambda )-g_{ij}^{1} (\lambda)\right|d\lambda  \le \delta_{i}^j, i,j=\overline{1,T}\right\}.$$
Here  $F_{1} ( \lambda ), G_1(\lambda)$ are known and fixed spectral densities, $W(\lambda)$ is an unknown spectral density, $p, \delta,\delta_{k}, p_k, k=\overline{1,T}$, $\delta_{i}^{j}, i,j=\overline{1,T}$, are given numbers, $P$ is a given positive-definite Hermitian matrices.

From the condition $0\in \partial \Delta _{D} (F^{0} ,G^{0} )$ we find the following equations which determine the least favourable spectral densities for these given sets of admissible spectral densities.

For the first pair $D_{\varepsilon}^1\times D_{1\delta}^1$ we have equations
\begin{equation} \label{eq_5_1}
(r_G^0(\lambda))^{*}(r_G^0(\lambda))^\top=(\alpha^{2} +\gamma_1(\lambda ))(F^{0} (\lambda
)+G^{0} (\lambda ))^{2},
\end{equation}
\begin{equation}\label{eq_5_2}
(r_F^0(\lambda))^{*}(r_F^0(\lambda))^\top=\beta^{2}\gamma_2(\lambda )(F^{0}(\lambda )+G^{0}(\lambda ))^{2},
\end{equation}
\begin{equation} \label{eq_5_3}
\frac{1}{2 \pi} \int_{-\pi}^{ \pi} \left|{\mathrm{Tr}}\, (G^0( \lambda )-G_{1}(\lambda )) \right|d\lambda =\delta,
\end{equation}
where $\alpha^{2}, \beta^{2}$ are Lagrange multipliers, $\gamma_1(\lambda )\le 0$ and $\gamma_1(\lambda )=0$ if ${\mathrm{Tr}}\,F^{0} (\lambda )>(1-\varepsilon ) {\mathrm{Tr}}\, F_{1} (\lambda )$, $ \left| \gamma_2( \lambda ) \right| \le 1$ and
\[\gamma_2( \lambda )={ \mathrm{sign}}\; ({\mathrm{Tr}}\, (G^{0} ( \lambda )-G_{1} ( \lambda ))): \; {\mathrm{Tr}}\, (G^{0} ( \lambda )-G_{1} ( \lambda )) \ne 0.\]

For the second pair $D_{\varepsilon}^2\times D_{1\delta}^2$ we have equations
\begin{equation} \label{eq_5_4}
(r_G^0(\lambda))^{*}(r_G^0(\lambda))^\top=(F^{0} (\lambda )+G^{0} (\lambda ))\left\{(\alpha_{k}^{2} +\gamma
_{k}^1 (\lambda ))\delta _{kl} \right\}_{k,l=1}^{T} (F^{0} (\lambda)+G^{0} (\lambda )),
\end{equation}
\begin{equation}\label{eq_5_5}
(r_F^0(\lambda))^{*}(r_F^0(\lambda))^\top=(F^{0}(\lambda )+G^{0} (\lambda ))\left\{\beta_{k}^{2} \gamma_{k}^2 (\lambda )\delta _{kl} \right\}_{k,l=1}^{T} (F^{0} (\lambda)+G^{0} (\lambda )),
\end{equation}
\begin{equation} \label{eq_5_6}
\frac{1}{2 \pi} \int_{- \pi}^{ \pi} \left|g^0_{kk} ( \lambda)-g_{kk}^{1} ( \lambda ) \right| d\lambda =\delta_{k},
\end{equation}
where  $\alpha _{k}^{2}, \beta_{k}^{2}$ are Lagrange multipliers, $\gamma_{k}^1(\lambda )\le 0$ and $\gamma_{k}^1 (\lambda )=0$ if $f_{kk}^{0}(\lambda )>(1-\varepsilon )f_{kk}^{1} (\lambda )$, $\left| \gamma^2_{k} ( \lambda ) \right| \le 1$ and
\[\gamma_{k}^2( \lambda )={ \mathrm{sign}}\;(g_{kk}^{0}( \lambda)-g_{kk}^{1} ( \lambda )): \; g_{kk}^{0} ( \lambda )-g_{kk}^{1}(\lambda ) \ne 0, \; k= \overline{1,T}.\]

For the third pair $D_{\varepsilon}^3\times D_{1\delta}^3$ we have equations
\begin{equation} \label{eq_5_7}
(r_G^0(\lambda))^{*}(r_G^0(\lambda))^\top=(\alpha^{2} +\gamma_1'(\lambda ))(F^{0} (\lambda )+G^{0} (\lambda ))B_{1}^{\top} (F^{0}(\lambda )+G^{0} (\lambda )),
\end{equation}
\begin{equation}\label{eq_5_8}
(r_F^0(\lambda))^{*}(r_F^0(\lambda))^\top=\beta^{2}\gamma_2'(\lambda)(F^{0}(\lambda )+G^{0}(\lambda))B_{2}^{\top}(F^{0}(\lambda)+G^{0}(\lambda)),
\end{equation}
\begin{equation} \label{eq_5_9}
\frac{1}{2 \pi} \int_{- \pi}^{ \pi} \left| \left \langle B_{2}, G^0( \lambda )-G_{1} ( \lambda ) \right \rangle \right|d\lambda
= \delta,
\end{equation}
where $\alpha^{2}, \beta^{2}$ are Lagrange multipliers, $\gamma_1' ( \lambda )\le 0$ and $\gamma_1' ( \lambda )=0$ if $\langle B_{1} ,F^{0} ( \lambda ) \rangle>(1- \varepsilon ) \langle B_{1} ,F_{1} ( \lambda ) \rangle$, $\left| \gamma_2' ( \lambda ) \right| \le 1$ and
\[\gamma_2' ( \lambda )={ \mathrm{sign}}\; \left \langle B_{2} ,G^{0} ( \lambda )-G_{1} ( \lambda ) \right \rangle : \; \left \langle B_{2} ,G^{0} ( \lambda )-G_{1} ( \lambda ) \right \rangle \ne 0.\]

For the fourth pair $D_{\varepsilon}^4\times D_{1\delta}^4$ we have equations
\begin{equation} \label{eq_5_10}
(r_G^0(\lambda))^{*}(r_G^0(\lambda))^\top=(F^{0} (\lambda )+G^{0}(\lambda ))(\vec{\alpha}\cdot \vec{\alpha}^{*}+\Gamma(\lambda))(F^{0} (\lambda )+G^{0} (\lambda )),
\end{equation}
\begin{equation}\label{eq_5_11}
(r_F^0(\lambda))^{*}(r_F^0(\lambda))^\top=(F^{0} (\lambda )+G^{0} (\lambda ))\left\{\beta_{ij}
\gamma_{ij}(\lambda ))\right\}_{i,j=1}^{T} (F^{0} (\lambda )+G^{0}(\lambda )),
\end{equation}
\begin{equation} \label{eq_5_12}
\frac{1}{2 \pi} \int_{- \pi}^{ \pi} \left|g^0_{ij}(\lambda)-g_{ij}^{1}( \lambda ) \right|d\lambda = \delta_{i}^{j},
\end{equation}
where  $\vec{\alpha}, \beta_{ij}$ are Lagrange multipliers, $\Gamma(\lambda )\le 0$ and $\Gamma(\lambda )=0$ if $F^{0}(\lambda )>(1-\varepsilon )F_{1} (\lambda )$, $\left| \gamma_{ij} ( \lambda ) \right| \le 1$ and
\[\gamma_{ij} ( \lambda )= \frac{g_{ij}^{0} ( \lambda )-g_{ij}^{1} (\lambda )}{ \left|g_{ij}^{0} ( \lambda )-g_{ij}^{1}(\lambda) \right|} : \; g_{ij}^{0} ( \lambda )-g_{ij}^{1} ( \lambda ) \ne 0, \; i,j= \overline{1,T}.\]

The following theorem and corollaries hold true.

\begin{thm}
Let the minimality condition (\ref{minimal-Unc}) hold true. The least favorable spectral densities  $F^0(\lambda)$, $G^0(\lambda)$  in the classes $D_{\varepsilon}\times D_{1\delta}$ for the optimal linear extrapolation of the functional $A\vec{\xi}$ are determined by relations
(\ref{eq_5_1}) -- (\ref{eq_5_3}) for the first pair $D_{\varepsilon}^1\times D_{1\delta}^1$ of sets of admissible spectral densities;
(\ref{eq_5_4}) -- (\ref{eq_5_6}) for the second pair $D_{\varepsilon}^2\times D_{1\delta}^2$ of sets of admissible spectral densities;
(\ref{eq_5_7}) -- (\ref{eq_5_9}) for the third pair $D_{\varepsilon}^3\times D_{1\delta}^3$ of sets of admissible spectral densities;
(\ref{eq_5_10}) -- (\ref{eq_5_12}) for the fourth pair $D_{\varepsilon}^4\times D_{1\delta}^4$ of sets of admissible spectral densities;
constrained optimization problem (\ref{extrem}) and restrictions  on densities from the corresponding classes $D_{\varepsilon}\times D_{1\delta}$.  The minimax-robust spectral characteristic of the optimal estimate of the functional $A\vec{\xi}$ is determined by the formula (\ref{sp-ch2}).
\end{thm}

\begin{nas}
Let the minimality condition (\ref{minimal1}) hold true. The least favorable spectral densities $F^{0}(\lambda)$ in the classes $D_\varepsilon^{k}$, $k=1,2,3,4$, for the optimal linear extrapolation of the functional  $A\vec{\xi}$ from observations of the sequence $\vec{\xi}(j)$ at points  $j\in\mathbb{Z}_{-}\backslash S $, where $S=\bigcup\limits_{l=1}^{s}\{-M_{l}-N_l, \ldots,-M_{l} \},$ are determined by the following  equations, respectively,
\begin{equation}
((C^{0}(\lambda) )^{\top} )^{*}\cdot(C^{0}(\lambda))^{\top}=(\alpha^{2}+\gamma(\lambda ))(F^{0} (\lambda ))^{2},
\end{equation}
\begin{equation}
((C^{0}(\lambda) )^{\top} )^{*}\cdot(C^{0}(\lambda) )^{\top}=F^{0} (\lambda )\left\{(\alpha_{k}^{2}+\gamma _{k} (\lambda ))\delta_{kl}\right\}_{k,l=1}^{T}F^{0} (\lambda ),
\end{equation}
\begin{equation}
((C^{0}(\lambda) )^{\top} )^{*}\cdot(C^{0}(\lambda) )^{\top}=(\alpha^{2}+\gamma' (\lambda)) F^{0} (\lambda )B_1^\top F^{0} (\lambda ),
\end{equation}
\begin{equation}
((C^{0}(\lambda) )^{\top} )^{*}\cdot(C^{0}(\lambda) )^{\top}=F^{0} (\lambda ) (\vec{\alpha}\cdot \vec{\alpha}^{*}+\Gamma(\lambda ))F^{0} (\lambda ),
\end{equation}
constrained optimization problem (\ref{extrem2}) and restrictions  on densities from the corresponding classes  $D_\varepsilon^{k}$, $k=1,2,3,4$. The minimax spectral characteristic of the optimal estimate of the functional $A\vec{\xi}$ is determined by the formula (\ref{spchar2}).
\end{nas}

\begin{nas}
Let the minimality condition (\ref{minimal1}) hold true. The least favorable spectral densities $F^{0}(\lambda)$ in the classes $D_{1\delta}^{k}$, $k=1,2,3,4$, for the optimal linear extrapolation of the functional  $A\vec{\xi}$ from observations of the sequence $\vec{\xi}(j)$ at points  $j\in\mathbb{Z}_{-}\backslash S $, where $S=\bigcup\limits_{l=1}^{s}\{-M_{l}-N_l, \ldots,-M_{l} \},$ are determined by the following  equations, respectively,
\begin{equation}
((C^{0}(\lambda) )^{\top} )^{*}\cdot(C^{0}(\lambda))^{\top}=\beta^{2} \gamma_2( \lambda )(F^{0} ( \lambda ))^{2},
\end{equation}
\begin{equation}
((C^{0}(\lambda) )^{\top} )^{*}\cdot(C^{0}(\lambda) )^{\top}=F^{0} ( \lambda ) \left \{ \beta_{k}^{2} \gamma^2_{k} ( \lambda ) \delta_{kl} \right \}_{k,l=1}^{T} F^{0} ( \lambda ),
\end{equation}
\begin{equation}
((C^{0}(\lambda) )^{\top} )^{*}\cdot(C^{0}(\lambda) )^{\top}=\beta^{2} \gamma_2'( \lambda )F^{0} ( \lambda )B_{2}^{ \top} F^{0} ( \lambda ),
\end{equation}
\begin{equation}
((C^{0}(\lambda) )^{\top} )^{*}\cdot(C^{0}(\lambda) )^{\top}=F^{0} ( \lambda ) \left \{ \beta_{ij}( \lambda ) \gamma_{ij} ( \lambda ) \right \}_{i,j=1}^{T} F^{0} ( \lambda ),
\end{equation}
constrained optimization problem (\ref{extrem2}) and the following restrictions on densities from the corresponding classes  $D_{1\delta}^{k}$, $k=1,2,3,4$, respectively,
\begin{equation}
\frac{1}{2 \pi} \int_{-\pi}^{ \pi} \left|{\mathrm{Tr}}\, (F^0( \lambda )-G_{1}(\lambda )) \right|d\lambda =\delta,
\end{equation}
\begin{equation}
\frac{1}{2 \pi} \int_{- \pi}^{ \pi} \left|f^0_{kk} ( \lambda)-g_{kk}^{1} ( \lambda ) \right| d\lambda =\delta_{k},
\end{equation}
\begin{equation}
\frac{1}{2 \pi} \int_{- \pi}^{ \pi} \left| \left \langle B_{2}, F^0( \lambda )-G_{1} ( \lambda ) \right \rangle \right|d\lambda
= \delta,
\end{equation}
\begin{equation}
\frac{1}{2 \pi} \int_{- \pi}^{ \pi} \left|f^0_{ij}(\lambda)-g_{ij}^{1}( \lambda ) \right|d\lambda = \delta_{i}^{j}.
\end{equation}
The minimax spectral characteristic of the optimal estimate of the functional $A\vec{\xi}$ is determined by the formula (\ref{spchar2}).
\end{nas}

\section{Conclusions}
In this article we describe methods of the mean-square optimal linear extrapolation of functionals which depend on the  unknown values of a multidimensional stationary sequence.
 Estimates are based on observations of the sequence with an additive stationary noise sequence.
 We develop methods of finding the optimal estimates of the functionals in the case of missing observations.
The problem is investigated in the case of spectral certainty, where the spectral densities of the sequences are exactly known.
In this case  we propose an approach based on the Hilbert space projection method.
We derive formulas for calculating the spectral characteristics and the mean-square errors of the estimates of the functionals.
In the case of spectral uncertainty, where the spectral densities of the sequences are not exactly known while sets of admissible spectral densities are given,
the minimax (robust) method of estimation is applied.
This method allows us to find estimates that minimize the maximum values of the mean-square errors
of estimates for all spectral density matrices from a given class of admissible spectral density matrices and derive
relations which determine the least favourable spectral density matrices.
These least favourable spectral density
matrices are solutions of the optimization problem
$\Delta_D(F,G)=-\Delta(h(F^0,G^0);F,G)+\delta((F,G)\left|D_F\times D_G\right.)\rightarrow \inf$.
These solutions are characterized by the condition $0 \in \partial\Delta_D(F^0,G^0)$, where $\partial\Delta_D(F^0,G^0)$ is the subdifferential of the
convex functional $\Delta_D(F,G)$ at point $(F^0,G^0)$.
The form of the functional $\Delta(h(F^0,G^0);F,G)$ is convenient for
application of the Lagrange method of indefinite multipliers for finding solution to the optimization problem. The
complexity of the problem is determined by the complexity of calculation of the subdifferential of the convex
functional $\Delta_D(F,G)$. Making use of the method of Lagrange multipliers and the form of subdifferentials of the
indicator functions we describe relations that determine the least favourable spectral densities in some special
classes of spectral densities.
These are: classes $D_0$ of densities with the moment restrictions, classes $D_{1\delta}$ which
describe the ``${\delta}$-neighborhood'' models in the space $L_1$ of a given bounded spectral density, classes
$D_U^V$ which
describe the "strip" models of a given bounded spectral density, classes $D_{\varepsilon}$ which describe the "${\varepsilon}$-contamination"
models of spectral densities.

\end{document}